\documentclass[12pt,twoside]{article}

\font\ninerm=cmr10

\newcommand{\ph}{\varphi}
\newcommand{\dis}{\displaystyle}
\newcommand{\bqr}{\begin{eqnarray}}
\newcommand{\bqre}{\begin{eqnarray*}}
\newcommand{\eqr}{\end{eqnarray}}
\newcommand{\eqre}{\end{eqnarray*}}

\title{Averaging lemmas with a force term in the transport equation}
\author{F. Berthelin and S. Junca}
\date{{\it \small Laboratoire J.A. Dieudonn\'e, CNRS UMR 6621, 
        \\ Universit\'e de Nice Sophia-Antipolis, 
       \\ Parc Valrose, 06108, Nice, France,
\\ bertheli@unice.fr, junca@unice.fr}}

\begin{document}

\font\bba=msbm10
\font\bbb=msbm8
%\font\bbb=msbm10 scaled 800
\font\bbc=msbm6
%\font\bbc=msbm10 scaled 600
\newfam\bbfam
\textfont\bbfam=\bba
\scriptfont\bbfam=\bbb
\scriptscriptfont\bbfam=\bbc
\def\bb{\fam\bbfam\bba}

\def\N{{\bb N}}
\def\Z{{\bb Z}}
\def\r{{\bb R}}
\def\C{{\bb C}}
\def\1{{1\hspace{-1.2mm}{\rm I}}}%\mbox{\normalsize I}}}
\def\emptyset{\hbox{$\displaystyle/\kern -5.97pt\circ$}}
\def\diam{\mbox{{\normalsize diam}}}
\def\supp{\mbox{{\normalsize supp}}}
\def\lip{\mbox{{\normalsize Lip}}}
\def\stackunder#1#2{\mathop{#1}\limits_{#2}}
\def\stackover#1#2{\mathop{#1}\limits^{#2}}
\def\CQFD{\unskip\kern 6pt\penalty 500%
\raise -1pt\hbox{\vrule\vbox to 8pt{\hrule width 6pt\vfill\hrule}\vrule}}
\def\ess{\mathop{\mbox{\rm ess}}}
\def\limsup{\mathop{\overline{\lim}}\limits}
\def\liminf{\mathop{\underline{\lim}}\limits}
\def\dv{\mathop{\mbox{\rm div}}\nolimits}
\def\sgn{\mathop{\mbox{\rm sgn}}}
\def\iint{\mathop{\int\mkern -12mu\int}}
\def\iiint{\mathop{\int\mkern -12mu\int\mkern -12mu\int}}
\def\iiiint{\mathop{\int\mkern -12mu\int\mkern -12mu\int\mkern -12mu\int}}
\def\tq{\mathrel{;\,}}
\def\pref#1{(\ref{#1})}
\def\pr{\noindent {\bf Proof. }}
\let\dsp=\displaystyle
\def\convol{\mathop{*}\limits}
\def\tr{\mathop{\mbox{\rm tr}}}
\def\meas{\mathop{\rm meas}}
\newcommand{\integr}[1]{\int_\r \int_\r #1  \,dx \,d\xi}
\newcommand{\integ}[1]{\int_\r  #1 \,d\xi}
\def\limd{\mathop{{\lim}}\limits}
\def\tod{\mathop{{\to}}\limits}
\def\supd{\mathop{{\sup}}\limits}
\def\mind{\mathop{{\min}}\limits}
\def\cupd{\mathop{{\cup}}\limits}
\def\simd{\mathop{{\sim}}\limits}
\def\ds{\displaystyle}
\def\eps{\varepsilon}

\newtheorem{Th}{Theorem}
\newtheorem{Prop}{Proposition}
\newtheorem{Lemma}{Lemma}
\newtheorem{Defin}{Definition}
\newtheorem{Cor}{Corollary}
\newtheorem{Rk}{\sl Remark}[section]
\newtheorem{Ex}{\sl Example}[section]

\def\theequation{\thesection.\arabic{equation}}
\def\thesection{\arabic{section}}
\def\thesubsection{\arabic{section}.\arabic{subsection}}
\def\thesubsubsection{\arabic{section}.\arabic{subsection}.\arabic{subsubsection}}
\newcommand\Section{%
\def\thesubsection{\arabic{section}}  
\setcounter{Th}{0}
\setcounter{Rk}{0}
\setcounter{Ex}{0}
\setcounter{equation}{0}\section}
\newcommand\Subsection{%
\def\thesubsection{\arabic{section}.\arabic{subsection}}
\subsection}
\newcommand\Subsubsection{\subsubsection}

\newcommand{\alali}{$\mbox{ }$\\} 
\maketitle

 %\begin{center}
	
 %\end{center}
 %\vspace{1cm}

%\def\abstractname{R\'esum\'e}
\begin{abstract}
We obtain several averaging lemmas for transport operator with a force term.
These lemmas  improve the regularity yet known 
by  not considering the 
force term as part of an arbitrary right-hand side.
Two methods are used: local variable changes or stationary phase.
These new results are subjected to two non degeneracy assumptions.
%, quantified by two parameters, the classical condition for the velocity field  and a new condition involving  the velocity field together with the force field. 
We characterize the optimal conditions of these assumptions
 to compare the obtained regularities 
according to the space and velocity variables.
Our results are mainly in  $L^2$,  and  for constant force, in $L^p$ for $1<p \leq 2$.
%We give several results using various tools :
%Fourier Series and Transform, Stationnary Phase Method, Hardy Space.
\end{abstract}
\vspace{4mm}

\def\abstractname{R\'esum\'e}
\begin{abstract}
Nous obtenons plusieurs lemmes de moyenne pour des \'equations de transport avec un terme de force.
Ces r\'esultats am\'eliorent la r\'egularit\'e connue en ne consid\'erant pas le terme de force
comme un terme source arbitraire.
Deux techniques sont utilis\'ees : des changements de variables locaux ou des phases stationnaires.
Ces r\'esultats sont quantifi\'ees par deux hypoth\`eses de non d\'eg\'en\'erescence.
Nous caract\'erisons les conditions optimales de ces hypoth\`eses 
pour comparer les r\'egularit\'es obtenues, par rapport aux variables d'espace et de vitesse.
Les r\'esultats sont principalement dans 
 $L^2$,  et pour le cas constant, dans $L^p$ pour $1<p \leq 2$.
%We give several results using various tools :
%Fourier Series and Transform, Stationnary Phase Method, Hardy Space.
\end{abstract}

\medskip

%35L65, 76N15, 35B35, 82C40
\vspace{1cm}

%\noindent {\it Work partially supported by 
%European TMR network projects\\
%Hyperbolic Conservation Laws contract \# ERBFMRXCT960033\\
%Kinetic theory contract \# ERBFMRXCT970157
%}

 \newpage
\baselineskip=12pt
\parskip=0pt plus 1pt

\tableofcontents

\Section{Introduction}

Averaging lemma is a major tool to get compactness from a kinetic equation.
(\cite{DPL1}, ...).
Such results have been used in a lot of papers during these last years.
Among this literature, an important result using
an averaging lemma as
a key argument is the proof of the hydrodynamic limits of the 
Boltzmann or BGK equations to the incompressible Euler or Navier-Stokes equations (\cite{GSR}).
 Another major application consists in obtaining the compactness
 for nonlinear scalar conservation  laws  (in \cite{LPT94}) 
which allows, for instance,
 to study the propagation of high frequency waves (\cite{CJR}). 

\noindent Basically, averaging lemma is a result which says that 
the macroscopic quantities $\ds \int f(t,x,v) \psi(v) \,dv$
have a better regularity with respect to  $(t,x)$ than 
the microscopic quantity $f(t,x,v)$ where $f$ is solution of a kinetic equation.

\noindent For example, in \cite{DPLM}  and \cite{Bez}, the
following result is established.
\medskip 
\\
% \begin{Thn}
{\bf Theorem [DiPerna, Lions, Meyer -- B\'ezard]} %     \label{thref}
\alali
{\it 
Let $f$, $g_k \in L^p(\r_t\times
\r^N_x\times\r^M_v)$ with $1<p \leq 2$ such that
\begin{eqnarray} \label{avecderivee}
\partial_t f+\dv_x [a(v) f]= \sum_{|k|\leq m}
\partial_v^k g_k,
\end{eqnarray}
with $a \in W^{m,\infty}(\r^M,\r^N)$ 
for $m\in \N$.
% $a$
%non d\'eg\'en\'er\'ee.
Let $\psi \in W^{m,\infty}(\r^M)$  with compact support.
Let $A>0$ such that the support of $\psi $ is included in $[-A,A]^M$.
We assume the following non-degeneracy for $a(.)$:
there exists $0 < \alpha \leq 1$ and $C>0$ such that
for any $(u,\sigma) \in S^{N}$ and $\eps>0$,
\begin{eqnarray*} 
       \nonumber
\textrm{meas } \left(\{v \in [-A,A]^M %\textrm{Supp} \psi 
    \, ; \, u-\eps<a(v)\cdot \sigma <u+\eps\} \right)\leq
C\eps^\alpha.
\end{eqnarray*}
Then $$ \ds \rho_\psi(t,x)=\int_{\r^M} f(t,x,v) \psi(v) \,dv$$ is in $W^{s,p}(\r_t\times\r^N_x)$ where
$s= \frac{\alpha}{(m+1)p'}$, $p'$ being the conjugated exponent for $p$.}\\
%We have besides}
%\begin{equation}
%\|(1-\Delta_{t,x})^{s/2} \rho_\psi \|_{L^p(\r_t\times\r^N_x)} \leq 
%C \left( \|f\|_{L^p(\r^N_x\times\r^M_v)}+
%\sum_\alpha \|g_\alpha\|_{L^p(\r^N_x\times\r^M_v)}\right).
%\end{equation}
%} % \end{Thn}

\noindent Regarding equation (\ref{avecderivee}),
the obtained regularity is proved to be optimal, see \cite{Lio3} and \cite{Lio4}.
In \cite{Ger2}, the gain of a half-derivative in $L^2$ context was proved
as optimal.
A study in the case of a full derivative with respect to $x$
in the second member is done in \cite{PS}.
We also refer to \cite{Ger1} and \cite{BD} for other results about averaging lemmas.
Regularity of $f$ itself is also challenging, for example by assuming some regularity
with respect to $v$, see \cite{Bou}, \cite{JP} and \cite{BB5} for such results.\\

\noindent Theorem here above says for example with $m=1$ that for the equation
\begin{equation}  \label{free}
\partial_t f + a(v) \cdot \nabla_x f  =g - F(t,x,v) \cdot \nabla_v \tilde g,
\end{equation}
the obtained regularity is $W^{s,p}(\r_t\times\r^N_x)$ with
$s= \frac{\alpha}{2p'}$.
When we consider equation 
\begin{equation} \label{generalforcesecond}
\partial_t f+ a(v) \cdot \nabla_x f +F(t,x,v) \cdot \nabla_v f=g,
\end{equation}
that is to say that $\tilde g=f$, it is classical to consider the term
$F(t,x,v) \cdot \nabla_v f$ being part of the right-hand side and to obtain 
the regularity $W^{s,p}(\r_t\times\r^N_x)$ with
$s= \frac{\alpha}{2p'}$.
But for (\ref{generalforcesecond}), the derivation with respect to $v$ is
only on $f$ through the transport equation and not on an arbitrary term $\tilde g$.
That is to say, the  conventional method is losing information because this term is part of characteristics
and the right-hand side terms are in $L^2$, % that is to say 
i.e. for $m=0$, and the obtained regularity should be 
$W^{s,p}(\r_t\times\r^N_x)$ with
$s= \frac{\alpha}{p'}$.\\
This is the first motivation of this paper 
and one of the result we get.\\
Few other papers deal with averaging lemma avoiding
to consider the acceleration term as a source, 
namely  \cite{GG}, \cite{Gol}. 
But they are based on a transversality assumption on $a(.)$
restricting the generality to the case $\alpha=1$.

% proven in Theorem \ref{lemmemoyenne} below.

\noindent Notations for (\ref{generalforcesecond}) are
$f(t,x,v) \in \r$ with $t \in \r$, $x \in \r^N$, $v \in \r^M$,
$a : \r^M \to \r^N$, $F: \r \times \r^N \times \r^M \to \r^M$ and
$$a(v) \cdot \nabla_x f=\sum_{i=1}^N a_i(v)\, \partial_{x_i} f,\quad
F(t,x,v) \cdot \nabla_v f = \sum_{i=1}^M F_i(t,x,v)\, \partial_{v_i} f.$$

\noindent In this paper, we will prove the following averaging lemmas on equation
(\ref{generalforcesecond}).
\begin{Th}[$L^2$ result] \label{lemmemoyenne}
%--------------------------
$\mbox{ }$\\
Let $a \in C^{N+3}(\r^M_v,\r^N_x)$, $F \in C^{N+3}(\r_t\times\r^N_x\times\r_v^M,\r^M_v)$,
$f, g \in L^2(\r_t\times\r^N_x\times\r_v^M)$, 
satisfying
(\ref{generalforcesecond}). Let  $A > 0$ and $\ds \psi \in 
C^{N+2}_c(\r_v^M)$ be such that the support of $\psi $ is included in $[-A,A]^M$.
We assume that
there exists  $0<\alpha \leq 1$ and $C>0$ such that
for any $(u,\sigma) \in S^{N}$ and $\eps>0$,
\begin{equation} 
       \label{alpha}
\textrm{meas } \left(\{v \in [-A,A]^M %\textrm{Supp} \psi 
    \, ; \, u-\eps<a(v)\cdot \sigma <u+\eps\} \right)\leq
C\eps^\alpha.
\end{equation}
%(L^2\cap L^\infty)(\r_v^M, W_c^{N+1+\alpha/4+\eta,\infty}(\r^{N+1}_{tx})),$ 
%with $\eta>0$,
Then the  averaging 
$$
\rho_\psi(t,x)=\int_{\r^M} f(t,x,v) \psi(v)\,dv
$$
is in $H^{\alpha/2}_{loc}(\r_t\times\r^N_x)$.
\end{Th}

\begin{Rk}
 We notice that we obtain $\alpha/2$ instead of the well known $\alpha/4$ when the acceleration term $F\cdot \nabla_x f$ is considered as a right hand side with no particular relation to $f$.
\end{Rk}

\begin{Rk}
For Vlasov equation, the classical application of averaging lemma is
the 
DiPerna, Lions, Meyer Theorem 
which
 gives the compactness
for $\rho_\psi$ with an operator of the kind (\ref{generalforcesecond})
applying the result with
$g_1=-F \cdot f$ when $F \in L^\infty_{loc}$.
More precisely, if $f^n$, $g_0^n$ and $g_1^n=-F_n \cdot f^n$ are
solutions of  (\ref{avecderivee}) with some bounds in $L^p$, then
$\ds \rho_\psi^n \textrm{ is bounded in } W^{s,p}(\r_t\times\r^N_x) \textrm{ with }
s=\frac{\alpha}{2p'},$ and thus is compact in
$W^{s',p}(\r_t\times\r^N_x)$ with $s'<s$.
For $p=2$, it is compact in $H^{s'}(\r_t\times\r^N_x)$ with $s'<\frac{\alpha}{4}$.
By this way, paper \cite{DPL2} proves
the existence of weak solutions for Vlasov-Maxwell.
With  Theorem 1, the obtained compactness is in 
$H^{s'}_{loc}(\r_t\times\r^N_x)$ with $s'<\frac{\alpha}{2}$.
\end{Rk}

\noindent When the force is constant, we obtain a global regularity 
    result with a less smooth test function. 
   
\begin{Th}[$L^2$ result with $F$ constant]  \label{FODELPNS}
%--------------------------
$\mbox{ }$\\
Let $a  \in C^\gamma(\r^M_v,\r^N_x)$, $F(t,x,v)=F \in \r^M$, $F \neq 0$,
$f$, $g \in L^2(\r_t\times \r^N_x\times\r^{M}_v)$ satisfying
(\ref{generalforcesecond}) where
we assume that function $a(\cdot)$ satisfies the following condition with $\gamma$, which is 
 a positive integer, such that
\\ $ \dis %\bqre
 \quad \forall (v,\sigma) \in  \r^M  \times S^{N},      
     \quad  \sigma =(\sigma_0,\sigma_1,\cdots,\sigma_N), 
       \quad \widetilde{\sigma}= (\sigma_1,\cdots,\sigma_N),
 $% \eqre
\bqr \label{gamma}
  |\sigma_0  + a(v).\widetilde{\sigma}| 
         +  
 \displaystyle{ \sum_{k=1}^{\gamma -1}
  \left| (F\cdot \nabla_v)^k a(v)\cdot \widetilde{\sigma} \right| }
  > 0.
&\qquad (\gamma ND)
 \eqr
Let $\psi\in C^{1}_c(\r^M_v)$, 
%  \in H^1(\r^N_v,\r) \cap C^1(\r^N_v,\r)$
then %  $\rho_\psi$ belongs in $H^{1/\gamma}(\r^d_X)$.
the   averaging 
$$ 
\rho_\psi(t,x)=\int_{\r^M} f(t,x,v) \psi(v)\,dv
   $$ 
is in $H^{1/\gamma}(\r_t \times \r^N_x)$.
% and there exists
% $K>0$
%independant of $f$ and $g$ such that
%\begin{equation} \label{resultat1}
%\|\rho_\psi\|_{H^{1/2}(\r^+\times\r^N_x)}\leq K
% \left( \|f\|_{L^2(\r^+\times\r^N_x\times\r^N_v)}+
%\|g\|_{L^2(\r^+\times\r^N_x\times\r^N_v)} \right)
%^{1/2} \|f\|_{L^2(\r^+\times\r^N_x\times\r^N_v)}^{1/2}.
%\end{equation}
\end{Th}

%
%\medskip
%
\begin{Rk}
 The proof of  Theorem \ref{FODELPNS} is not  valid when $F =0$. So this theorem does not give an   averaging Lemma for the   kinetic equation $\partial_t f + a(v) \cdot \nabla_x f = g $. 
\end{Rk}

\begin{Rk}
The case of a nonzero constant force field is not without interest, as it appears for instance when considering 
gravity effects in the kinetic theory of neutral gases.
\end{Rk}

\begin{Rk} {\bf [$M=1$, one dimensional velocity ]} $\mbox{ }$ % \\
 \begin{enumerate}
\item
% 1)
 The Sobolev estimate for $\rho_\psi$ comes from optimal bounds 
 in  stationary  phase lemma. Then, with only $f,g \in L^2$ and $M=1$, 
 we expect  Theorem  \ref{FODELPNS} to
 give the best Sobolev's exponent.
 \item  % 2) 
     Since $\gamma \geq N+1$ (see Proposition \ref{mingamma} for this inequality), 
    with only $f,g \in L^2$, 
   we expect  $ \rho_\psi$ to belong at most  to $H^{1/(N+1)}(\r^{N+1}_X)$
 when $M=1$.
\item With scalar velocity, the condition $(\gamma ND)$ 
  is similar  to  a non degeneracy condition given in \cite{Gof} about 
averaging for operators with real principal symbols. 
 More precisely it is the condition (5) of Theorem 4 with $t=v$ and $\xi_0=F$ in
\cite{Gof}.
 But our result yields a better smoothing effect, the gain of regularity for the
average is  $1/\gamma$ 
 instead of  $1/(2(\gamma -1))$ in \cite{Gof}.
\end{enumerate}
 \end{Rk} 

\noindent Next Theorem is a comparison between
 the two previous results. It shows that Theorem \ref{lemmemoyenne}
 does not give the best Sobolev exponent when $M=1$
 and that Theorem \ref{FODELPNS} is not 
 optimal for $M>1$.

\begin{Th}  \label{Pcomparison}
 For $N \geq 2$ and $M=1$, Theorem  \ref{FODELPNS} gives a stronger smoothing 
    effect than Theorem \ref{lemmemoyenne}
 for the best  $\gamma= \gamma_{opt}$ compared with the best $\alpha=\alpha_{opt}$
  since
  \bqre 
       \frac{1}{ \gamma_{opt}} = \frac{1}{N+1}
   & > &
  \frac{\alpha_{opt}}{2} =  \frac{1}{2N}.
  \eqre
Conversely, for $N=M$, Theorem \ref{lemmemoyenne} can give one half derivative 
   with the best $\alpha=1$.
\end{Th}
\begin{Rk}
For scalar velocity ($v \in \r,$ $M=1$), we characterize in Theorem \ref{Pcomparison} the best parameter $\alpha$ for the classical non degeneracy condition, namely condition (\ref{alpha}).  This characterization is mentioned in few works, see \cite{LPT94,JaX09}, but the proof of  optimality is a new result. This kind of characterization also gives new results for scalar conservation laws, see \cite{Ju}.
\end{Rk}

\noindent Finally, we find out two results in $L^p$ framework.

\begin{Th}[First $L^p$ result with $F$ constant] 
%--------------------------
$\mbox{ }$\\
Let $a \in C^{N+3}(\r^M_v,\r^N_x)$, $F(t,x,v)=F \in \r^M_v$,
$f, g \in L^p(\r_t\times\r^N_x\times\r_v^M)$, 
satisfying
(\ref{generalforcesecond}). Let  $A > 0$ and $\ds \psi \in 
C^{N+2}_c(\r_v^M)$ be such that the support of $\psi $ is included in $[-A,A]^M$.
We assume that
there exists  $0<\alpha \leq 1$ and $C>0$ such that
for any $(u,\sigma) \in S^{N}$ and $\eps>0$,
\begin{equation} 
       \label{alpha2}
\textrm{meas } \left(\{v \in [-A,A]^M %\textrm{Supp} \psi 
    \, ; \, u-\eps<a(v)\cdot \sigma <u+\eps\} \right)\leq
C\eps^\alpha.
\end{equation}
%(L^2\cap L^\infty)(\r_v^M, W_c^{N+1+\alpha/4+\eta,\infty}(\r^{N+1}_{tx})),$ 
%with $\eta>0$,
Then the  averaging 
$$ 
\rho_\psi(t,x)=\int_{\r^M} f(t,x,v) \psi(v)\,dv
 $$
is in $W^{s,p}_{loc}(\r_t\times\r^N_x)$ with $\dis s=\frac{\alpha}{p'}$.
\end{Th}

\begin{Th}[Second $L^p$ result with $F$ constant]  \label{GLPS2}
%--------------------------
$\mbox{ }$\\
Let $a \in C^{\gamma}(\r^M_v,\r^N_x)$, $F(t,x,v)=F \in \r^M_v$, $F \neq 0$,
$f$, $g \in L^p(\r_t\times\r^N_x\times\r^M_v)$ 
with $1<p\leq 2$, satisfying 
(\ref{generalforcesecond}),
 where
we assume that $a(\cdot)$ satisfies the following condition with $\gamma$, which is 
 a positive integer, such that
$$\forall (v,\sigma) \in  \r^M  \times S^{N},      
     \quad  \sigma =(\sigma_0,\sigma_1,\cdots,\sigma_N), 
       \quad \widetilde{\sigma}= (\sigma_1,\cdots,\sigma_N), \quad $$
$$  |\sigma_0  + a(v).\widetilde{\sigma}| 
         +  
 \displaystyle{ \sum_{k=1}^{\gamma -1}
  \left| (F\cdot \nabla_v)^k a(v)\cdot \widetilde{\sigma} \right| }
  > 0.
\qquad (\gamma ND)$$
Let $\psi \in 
C^1_c(\r_v^M)$,
%, C_c^{N+1+\alpha/4+\eta}(\r^{N+1}_{tx}))$ with $\eta>0$.
then the averaging
$$
\rho_\psi(t,x)=\int_\r f(t,x,v) \psi(v)\,dv$$
is in  $W^{s,p}(\r_t\times\r^N_x)$ 
with $s=\frac{2}{\gamma p'}$.
\end{Th}

\begin{Rk}
These results are presented with time dependence because it is more useful 
in applications.\\
In the proof of next sections, we take the following notations.
We set $X=(t,x)$ and $b(v)=(1,a(v))$.
Then (\ref{generalforcesecond}) can be rewritten as follows:
\begin{equation}  \label{eqb}
 b(v) \cdot \nabla_X f +F(X,v) \cdot \nabla_v f=g,
\end{equation}
where  $ X \in \r^{N+1}$,  $v \in \r^M$. 
\end{Rk}

\noindent Here is how the paper is structured.\\
In Section 2, we prove Theorem 1 for a smooth force field. In Section 3, we prove Theorem 2 for a constant  and non zero force field.
In Section 4, we compare both results (Theorem 3)  and finally
in Section 5, we prove the extension to $L^p$ spaces for contant force (Theorem 4 and 5).

%%%%%%%%%%%%%%%%%%%%%%%%%%%%%%%%%%%%%%%%%%%%%%%%%%%%%%%%%%%%%%%%%%%%%%%%%%%%%%%
%%%%%%%%%%%%%%%%%%%%%%%%%%%%%%%%%%%%%%%%%%%%%%%%%%%%%%%%%%%%%%%%%%%%%%%%%%%%%%

\Section{First Theorem in the $L^2$ framework}

%%%%%%%%%%%%%%%%%%%%%%%%%%%%%%%%%%%%%%%%%%%%%%%%%%%%%%%%%%%%%%%%%%%%%%%%%%%%%%%
%%%%%%%%%%%%%%%%%%%%%%%%%%%%%%%%%%%%%%%%%%%%%%%%%%%%%%%%%%%%%%%%%%%%%%%%%%%%%%%
We first recall the following classical averaging lemma (see \cite{GLPS}, \cite{BGP}).
\begin{Prop}[Golse, Lions, Perthame, Sentis] \label{clas}
\alali
Let $a \in L^\infty_{loc}(\r^M,\r^N)$, $f, g \in L^2(\r_t\times\r^N_x\times\r_v^M)$, 
 such that
\begin{equation}  
\partial_t f + a(v) \cdot \nabla_x f  = g.
\end{equation}
Let $\psi \in L^\infty(\r^M_v)$, with compact support in some $[-A,A]^M$, such that
there exists $0<\alpha \leq 1$ and $C>0$ such that
\begin{equation} 
\textrm{meas } \left(\{v \in [-A,A]^M\, ; \, u-\eps<a(v)\cdot \sigma <u+\eps\} \right)\leq
C\eps^\alpha
\end{equation}
for any $(u,\sigma) \in S^{N}$ and $\eps >0$.
Then the averaging
$$
\rho_\psi(t,x)=\int_{\r^M} f(t,x,v) \psi(v)\,dv$$
is in $H^{\alpha/2}(\r_t\times\r^N_x)$
with the estimate 
$$\|\rho_\psi\|_{H^{\alpha/2}} \leq \tilde C(N) \left( \|\psi\|_{L^2} + \sqrt{K}\|\psi\|_{L^\infty}
 \right) \left( \|f\|_{L^2} + \|g\|_{L^2} \right).$$
\end{Prop}
We use this averaging lemma to prove an other result,
which deals with 
test function depending on $(t,x,v)$.
\begin{Prop}[Averaging Lemma with test function  in $(X,v)$] \label{LalXv} 
\alali
Let $a \in L^\infty_{loc}(\r^M_v,\r^N_x)$, $f, g \in L^2(\r_t\times\r^N_x\times\r_v^M)$, 
 such that
\begin{equation}  
\partial_t f + a(v) \cdot \nabla_x f  = g.
\end{equation}
Let $\ds \psi \in 
L^\infty_c(\r_v^M,W^{N+2,\infty}(\r_{tx}^{N+1}))$
with compact support with respect to $v$ in some $[-A,A]^M$.
We assume that
there exists  $0<\alpha \leq 1$ and $C>0$ such that
\begin{equation} 
\textrm{meas }\left(\{v\in [-A,A]^M\, ; \, u-\eps<a(v)\cdot \sigma <u+\eps\} \right)\leq
C\eps^\alpha
\end{equation}
for any $(u,\sigma) \in S^{N}$ and $\eps>0$.\\
Then, for any compact $K$, 
there exists a constant $C(N,K)$ such that the averaging
$$
\rho_\psi(t,x)=\int_\r f(t,x,v) \psi(t,x,v)\,dv$$
is in $H^{\alpha/2}_{loc}(\r_t\times\r^N_x)$
with the bound 
$$\|\rho_\psi\|_{H^{\alpha/2}_K}  \leq 
 C(N,K) \left( \|f\|_{L^2} + \|g\|_{L^2} \right) \| \psi \|_{
(L^2\cap L^\infty)_v(W^{N+2,\infty}_{tx})}.$$
\end{Prop}
\pr
We fix a compact $K$ on $X$. We take $\tilde{K}=[-S,S]^{N+1}$
such that $K \subset \tilde K$ and $\chi$ a $C^\infty$ function such that
$\chi =1$ on $K$ and $0$ outside $\tilde{K}$.
Finally, we set $\tilde{\psi} = \psi \chi$.\\
Since $\tilde{\psi}$ has a compact support with respect to $X$, 
we can extend it by periodicity in these variables.
Then the Fourier expansion with respect to $X$ gives
$$\tilde{\psi}(X,v)=\sum_{\beta \in \Z^{N+1}} c_\beta(v) e^{i S \beta \cdot X}.$$
We write this formula through
$$\tilde{\psi}(X,v)=\sum_{\beta \in \Z^{N+1}} \bigg((1+|\beta|^r) c_\beta(v)\bigg) \, \cdot \,
\frac{ e^{i S \beta \cdot X}}{1+|\beta|^r},$$
with $r=N/2+1$.
We set
$$\phi_\beta(X) = \frac{ e^{i S \beta \cdot X}}{1+|\beta|^r},  \, \, 
\textrm{ and } \, \, \psi_\beta(v) = (1+ |\beta|^r) c_\beta(v).$$
%The conditions $\phi_\beta \in C^1(\r\times\r^N)$, 
%$\phi_\beta$, $\phi_\beta' \in L^\infty (\r\times\r^N)$ are clear.
%The condition  $\ds \sum_{\beta \in \Z^{N+1}} |\phi_\beta(t,x)|^2$ converge
%is satisfied when $2r >N+1$.
%Now $$|\psi_\beta(v)| \leq C_5 |\beta|^r \sup_{t,x} |\psi(t,x,v)|$$
%and $\psi \in C^0_c(\r_t\times\r^N_x\times\r_v)$, 
%leads to $\psi_\beta \in L^2(\r)\cap L^\infty(\r)$.
%for the last condition, 
We use the decreasing of 
Fourier coefficients for $W^{N+2,\infty}(\r^{N+1}_X)$ function, that is to say that
$$
|c_\beta(v)| \leq  \frac{C_1}{(S|\beta|)^{N+2}}
\|\tilde\psi(\cdot,v)\|_{W_X^{N+2,\infty}}.
%\\
%& \leq & C \frac{1}{|\beta|^{N+3}} \| \psi\|_{L^\infty_v(W^{N+3,1}_{tx})},
%\\
%& \leq & C \frac{1}{|\beta|^{N+1+2r}} \| \psi\|_{L^\infty_v(W^{2N+4,1}_{tx})},
$$
%(We notice that this result for $W^{N+1,\infty}(\r)$ function is obtained by integration by parts
%and for the H\"older part $W^{\alpha/4+\eta,\infty}(\r)$,
%it comes from an estimate of this kind:
%\begin{eqnarray*}
%2 \left|\int_0^1 \psi(x) e^{iSnx} \,dx \right| & = & 
%\left|\int_0^1 (\psi(x)-\psi(x+\pi/(Sn)) e^{iSnx} \,dx \right|\\
%& = & \left( \frac{\pi}{Sn}\right)^{\alpha/4+\eta} C
%\end{eqnarray*}
%for function in one dimension.)
Thus we have
\begin{eqnarray}
&& \int_{\r^M} \sum_{\beta \in (\Z^{N+1})^*} |\psi_\beta(v)|^2 \,dv \nonumber \\
& \leq &  \int_{\r^M} \sum_{\beta \in (\Z^{N+1})^*} (1+|\beta|^{r})^2 |c_\beta(v)|^2 \,dv \nonumber\\
& \leq &  \frac{C_2}{S^{2N+4}} \int_{\r^M} \sum_{\beta \in (\Z^{N+1})^*}
 \frac{4|\beta|^{2r} }{|\beta|^{2N+4}} 
\|\tilde\psi(\cdot,v)\|_{W_X^{N+2,\infty}}^2
\,dv  \nonumber\\
%& \leq &  \frac{4C_2}{S^{N+1+\alpha/4+\eta}}\int_{\r^M} \sum_{\beta \in (\Z^{N+1})^*} 
%\frac{1}{|\beta|^{N+1+\alpha/2}} 
%\left(\int_{\r^{N+1}} |\partial^{N+1+\alpha/4+\eta}_{X_{j_\beta}...{X_{k_\beta}}} \psi(Y,v)| \,dY\right)^2
%\,dv \nonumber\\
& \leq &  \frac{4C_2}{S^{2N+4}}\sum_{\beta \in (\Z^{N+1})^*} \frac{1}{|\beta|^{N+2}} 
\| \psi \|_{L^2_v(W^{N+2,\infty}_{X})}^2  <+\infty.   \label{borne1}
\end{eqnarray}
%with $R \in L^2(\r)$ and since $v \mapsto \psi(t,x,v) \in C^{2r+N+1}(\r)$.
%We start by the case $\ds \psi(t,x,v)=\sum_{n=0}^{+\infty} \phi_n(t,x) \psi_n(v)$
%with $\phi_n \in C^1(\r\times\r^{N+1})$, 
%$\psi_n \in L^2(\r)\cap L^\infty(\r)$ and such that 
%$\phi_n$, $\phi_n' \in L^\infty (\r\times\r^{N+1})$,
%$\ds \sum_{n=0}^{+\infty} |\phi_n(t,x)|^2$ converge and 
%$\ds \int_\r \sum_{n=0}^{+\infty} |\psi_n(v)|^2 \,dv <+\infty$.
%$\ds \sum_{n=0}^{+\infty} \phi_n(t,x) \psi_n(v)$
%and $\ds \sum_{n=0}^{+\infty} \phi_n'(t,x) \psi_n'(v)$ converge uniformly on $\r\times\r^{N+1}\times\r$.
On $K$, we notice that
\begin{eqnarray*}
\rho_\psi(X)  &=&
\int_\r f(X,v) \psi(X,v) \,dv \, \chi(X)  \\
&=& \int_\r f(X,v) \tilde \psi(X,v) \,dv, \\
&=& \int_{\r^M} f(X,v) \sum_{\beta\in \Z^{N+1}} \phi_\beta(X) \psi_\beta(v)\,dv.
\end{eqnarray*}
To apply Fubini's Theorem, we need that,  for $a.e.$ X,
$$\int_{\r^M} \sum_{\beta\in (\Z^{N+1})^*} |f(X,v) \phi_\beta(X) \psi_\beta(v)| \,dv < +\infty.$$
It comes from
\begin{eqnarray*}
&&\int_{\r^M} \sum_{\beta\in (\Z^{N+1})^*} |f(X,v) \phi_\beta(X) \psi_\beta(v)| \,dv\\
&\leq& \int_{\r^M} |f(X,v)| \sum_{\beta\in (\Z^{N+1})^*} |\phi_\beta(X) \psi_\beta(v)| \,dv \\
&\leq & \sqrt{\int_{\r^M} |f(X,v)|^2 \,dv}
\sqrt{\int_{\r^M} \left( \sum_{\beta\in (\Z^{N+1})^*} |\phi_\beta(X) \psi_\beta(v)| \right)^2 \,dv}\\
&\leq & \|f(X,\cdot)\|_{L^2_{v}}
\sqrt{ \sum_{\beta\in (\Z^{N+1})^*} |\phi_\beta(X)|^2  \int_{\r^M}
\sum_{\beta\in (\Z^{N+1})^*} |\psi_\beta(v)|^2 \,dv}\\
&\leq & \|f(X,.\cdot)\|_{L^2_{v}}
\sqrt{ \sum_{\beta\in (\Z^{N+1})^*} \frac{1}{(1+|\beta|^{r})^2}  \int_{\r^M}
\sum_{\beta\in (\Z^{N+1})^*} |\psi_\beta(v)|^2 \,dv} \, <+\infty
\end{eqnarray*}
since $2r >N+1$ and from (\ref{borne1}).
Thus we can write, on $K$,
\begin{eqnarray*} \\
  \rho_\psi(X) &=&\sum_{\beta\in\Z^{N+1}} \phi_\beta(X) \rho_{\psi_\beta}(X),
 \\ \dis \rho_{\psi_\beta}(X)& = &\int_\r f(X,v) \psi_\beta(v)\,dv. 
\end{eqnarray*}
The classical averaging lemma (Proposition \ref{clas}) gives that
$$\|\rho_{\psi_\beta}\|_{H^{\alpha/2}_K} \leq \tilde C(N) \left( \|\psi_\beta\|_{L^2} + \sqrt{C}
\|\psi_\beta\|_{L^\infty} \right) \left( \|f\|_{L^2} + \|g\|_{L^2} \right).$$
We now use the following property:  
%$$\|u_1u_2\|_{H^s} \leq  C_2 \|u_1\|_{C^1}\|u_2\|_{L^2}+\|u_1\|_{L^\infty}\|u_2\|_{H^s}
%\leq C_3 \|u_1\|_{C^1}\|u_2\|_{H^s},$$
For $u_1 \in C^{s}(\Omega)$,  $u_2 \in H^s(\Omega)$, 
with $s \in ]0,1[$, 
with $\Omega$ a bounded open set of $\r^{N+1}$,
we have $u_1u_2 \in H^s(\Omega)$ with
$$\|u_1u_2\|_{H^s} \leq  C_3 \|u_1\|_{C^{s}}\|u_2\|_{H^s}.$$
This result gives, for $s=\alpha/2$,
\begin{eqnarray*}
&& \|\rho_\psi\|_{H^{\alpha/2}_K}\\ & \leq &
C_3 \sum_{\beta\in\Z^{N+1}} 
\|\phi_\beta\|_{C^{\alpha/2}_K}\|\rho_{\psi_\beta}\|_{H^{\alpha/2}_K}\\
& \leq &
C_4 \sum_{\beta\in(\Z^{N+1})^*} 
\|\phi_\beta\|_{C^{\alpha/2}_K} \left( \|\psi_\beta\|_{L^2} + \|\psi_\beta\|_{L^\infty} \right)
\left( \|f\|_{L^2} + \|g\|_{L^2} \right) + C_3 \|\rho_{\psi_0}\|_{H^{\alpha/2}_K}\\
& \leq & C_5 \left( \sum_{\beta\in(\Z^{N+1})^*} 
\frac{1}{|\beta|^{r-\alpha/2}} \frac{
\| \psi \|_{(L^2\cap L^\infty)_v(W^{N+2,\infty}_{X})}
}{|\beta|^{N+2-r} } 
\left( \|f\|_{L^2} + \|g\|_{L^2} \right)+ \|\tilde\psi\|_{C^1_K}\right)\\
& \leq & C_5 \left(\sum_{\beta\in(\Z^{N+1})^*} 
\frac{1}{|\beta|^{N+2-\alpha/2}}  
\left( \|f\|_{L^2} + \|g\|_{L^2} \right) 
\| \psi \|_{(L^2\cap L^\infty)_v(W^{N+2,\infty}_{X})} + \|\psi\|_{C^{N+2}_c}\right).
\end{eqnarray*}
Since $N+2-\alpha/2 >N+1$, the proof is completed. $\CQFD$\\

\noindent With this Proposition now stated, we can go into the proof of our first Theorem.\\

\noindent {\bf Proof of Theorem \ref{lemmemoyenne}.}
Let $K$ be a compact in $\r^{N+1}_X$. We set
${\cal K }= K \times [-A,A]^M$.
We perform locally a change in variables in order to rewrite equation (\ref{eqb})
without the term $\nabla_v f$ and to apply previous result.
For any $(X_0,v_0) \in {\cal K}$, 
using the characteristics since $b(v)=(1,a(v)) \neq 0$,
there exists ${\cal B}_{Xv} \subset {\cal K}$ a neighborhood of $(X_0,v_0)$
and a $C^{N+3}$ diffeomorphism 
$$ \begin{array}{cccl}
        \Phi_0\; :& {\cal B}_0  & \rightarrow  & {\cal B}^0, \\ 
                 & (X,w) & \mapsto  & \Phi_0(X,w)=(X,V_0(X,w)),
  \end{array}$$
such that on  $ {\cal B}_0$ we have
 \begin{eqnarray} \label{eqV0}
  b(V_0(X,w)) \cdot \nabla_X V_0(X,w) &=& F(X,V_0(X,w)).
 \end{eqnarray}
Let us explain more precisely how to define the diffeomorphism $\Phi_0$ 
from equation (\ref{eqV0}).
 Since $b(v)=(1,a(v))$, $X=(t,x)$ and $X_0=(t_0,x_0)$, 
 equation  (\ref{eqV0}) can be reformulated as a nonlinear hyperbolic system
 (where $w$ is a parameter)
 \begin{eqnarray} \label{eqV0t}
  \partial_t V_0(t,x;w) + a(V_0(t,x;w)) \cdot  \nabla_x V_0(t,x;w) &=& F(t,x,V_0(t,x;w)),
 \end{eqnarray}
 completed by the initial data
 \begin{eqnarray} \label{eqV0t0}
  V_0(t_0,x;w))  & = &   w.
 \end{eqnarray}
By the classical method of characteristics, for each $w$, there exists a neighborhood of $(t_0,x_0)$  where $V_0$ is well defined and smooth. The characteristics are smooth with respect to the parameter $w$, thus $V_0(t,x,w)$ is well defined on a neighborhood of $(t_0,x_0;v_0)$.    
 Notice that $\partial_w V_0(t_0,x;w)= id_{\r^M}$, with $id_{\r^M}$ the identity operator on $\r^M_v$, and $\det (D \Phi_0)= \det (\partial_w V_0)$,  so reducing if necessary the previous neighborhood,  $\Phi_0$ is a diffeomorphism on ${\cal B}_0$. 
\\
Denoting by $\tilde{f}_0(X,w)=f(X,V_0(X,w))$,  
$\tilde{g}_0(X,w)=g(X,V_0(X,w))$, 
$\tilde{b}_0(w)= b(V_0(X,w))$,
the equation (\ref{eqb}) rewrites
\begin{equation}
\tilde{b}_0(w) \cdot \nabla_X \tilde{f}_0 = \tilde{g}_0.
\end{equation}

\noindent Now, there exists a finite number of 
${\cal B}^{l}$ to recover this compact, i.e. 
there exists $\{(X_l,v_l)\}_{l=1,\cdots,L}$, with the associated diffeomorphim
 $\Phi_l: {\cal B}_l \rightarrow {\cal B}^l$, $\Phi_l(X,w)=(X,V_l(X,w))$, 
such that
${\cal K} \subset \cupd_{l=1,\cdots,L} {\cal B}^{l}$.
For this recovering, we use a partition of unity, we have 
\begin{eqnarray*} 
f(X,v) &=& f(X,v) \1_{\cal K}(X,v) =\sum_{l=1}^L f(X,v) \chi_l(X,v)  
\end{eqnarray*}
where  function $\chi_l$ are $C^\infty$ 
and have a compact support in ${\cal B}^{l}$.
\\
Denoting  again  by 
$\tilde{f}_l(X,w)=f(X,V_l(X,w))$,
$ \tilde{g}_l(X,w)=g(X,V_l(X,w))$,
$ \tilde{b}_l(w)= b(V_l(X,w))$ on ${\cal B}_l$
 and \begin{eqnarray*}
       I^l[X] &=& \{v \in \r^M \mbox{ such that } (X,v) \in B^l  \}, \\
       I_l[X]& =& \{w \in \r^M \mbox{ such that } (X,w) \in B_l  \},
    \end{eqnarray*}
 we have the following decomposition.
It is
\begin{eqnarray*}
 \rho_\psi(X) &=& \sum_{l=1}^L \int_{\r^M} f_l(X,v)\chi_l(X,v)  \psi(v) \,dv\\
&=& \sum_{l=1}^L \int_{I^l[X]} f(X,v) \chi_l(X,v) \psi(v) \,dv\\
& =& \sum_{l=1}^L \int_{I_l[X]} \tilde{f}(X,w)  \chi_l(X,V_l(X,w))
\psi(V_l(X,w)) J_l(X,w)\,dw.
\end{eqnarray*}
 where we can 
perform the variable change $v \mapsto w=V(X,v)$ on every neighborhood 
${\cal B}^l$ 
corresponding to $l$ and  denoting by $J_l(X,w)$ the associated jacobian, 
 i.e. $J_l = |\det D \Phi_l| = |\det \partial_w V_l|$.\\
We set 
$\overline{\psi_l}(X,w)=\chi_l(X,V(X,w))\psi(V(X,w)) J_l(X,w).$
Since $a$ and $F$ have $C^{N+3}$ regularity, $J_l$ has $C^{N+2}$ one.
Furthermore $\psi \in C^{N+2}_c$, thus $\overline{\psi_l} \in 
(L^2\cap L^\infty)_c(\r_v^M,W^{N+2,\infty}(\r^{N+1}_{X}))$.
We apply previous result, namely Proposition \ref{LalXv}, on  the averaging
$$
\rho_{\overline{\psi_l}}(X)=\int_{\r^M} \tilde f(X,w) \overline{\psi_l}(X,w)\,dw
\qquad \mbox{ which  is in } H^{\alpha/2}_{loc}(\r^{N+1}_X).$$
Finally  the inequality $\displaystyle
 \|\rho_\psi\|_{H^{\alpha/2}_K}  \leq  \sum_{l=1}^L \|\rho_{\overline{\psi_l}}\|_{H^{\alpha/2}_K}$
concludes the proof.\CQFD

%°°°°°°°°°°°°°°°°°°°°°°°°°°°°°°°°°°°°°°°°°°°°°°°°°°°°°°°°°°°°°°°°°°°°°°°°°°°°°
%°°°°°°°°°°°°°°°°°°°°°°°°°°°°°°°°°°°°°°°°°°°°°°°°°°°°°°°°°°°°°°°°°°°°°°°°°°°°°
%°°°°°°°°°°°°°°°°°°°°°°°°°°°°°°°°°°°°°°°°°°°°°°°°°°°°°°°°°°°°°°°°°°°°°°°°°°°°°
%°°°°°°°°°°°°°°°°°°°°°°°°°°°°°°°°°°°°°°°°°°°°°°°°°°°°°°°°°°°°°°°°°°°°°°°°°°°°°
%°°°°°°°°°°°°°°°°°°°°°°°°°°°°°°°°°°°°°°°°°°°°°°°°°°°°°°°°°°°°°°°°°°°°°°°°°°°°°
%°°°°°°°°°°°°°°°°°°°°°°°°°°°°°°°°°°°°°°°°°°°°°°°°°°°°°°°°°°°°°°°°°°°°°°°°°°°°°

\Section{Case of a constant force field} \label{sF}

%°°°°°°°°°°°°°°°°°°°°°°°°°°°°°°°°°°°°°°°°°°°°°°°°°°°°°°°°°°°°°°°°°°°°°°°°°°°°°
%°°°°°°°°°°°°°°°°°°°°°°°°°°°°°°°°°°°°°°°°°°°°°°°°°°°°°°°°°°°°°°°°°°°°°°°°°°°°°
%°°°°°°°°°°°°°°°°°°°°°°°°°°°°°°°°°°°°°°°°°°°°°°°°°°°°°°°°°°°°°°°°°°°°°°°°°°°°°
%°°°°°°°°°°°°°°°°°°°°°°°°°°°°°°°°°°°°°°°°°°°°°°°°°°°°°°°°°°°°°°°°°°°°°°°°°°°°°
%°°°°°°°°°°°°°°°°°°°°°°°°°°°°°°°°°°°°°°°°°°°°°°°°°°°°°°°°°°°°°°°°°°°°°°°°°°°°°
%°°°°°°°°°°°°°°°°°°°°°°°°°°°°°°°°°°°°°°°°°°°°°°°°°°°°°°°°°°°°°°°°°°°°°°°°°°°°°

When $F$ is a non zero  constant vector, we can obtain 
a different result.
 The way to get it
 is quite different and we have to be restricted to the case
of a constant force field.
  A key tool here is a generalized  uniform version of the 
  classical method of the stationary phase.
 %from  Elias, M. Stein, 
 % see \cite{Ste} for instance. 
  We work on equation (\ref{eqb}) with $F$ constant, $F \in \r^M$, $F \neq 0$.
 % \\
Let us denote a directional $v-$derivative along vector $F$ by
 \begin{eqnarray}
  D  & = &  F \cdot \nabla_v .
 \end{eqnarray}
The smoothing effect depends on   $(\gamma ND)$ assumption of 
Theorem  \ref{FODELPNS}.
 Indeed, it is exactly the following    non-degeneracy condition
 about $D$-derivatives of $b(.)$:
% $$
%\label{gND}
%\begin{array}{l}
%\forall \sigma \in S^{N},\mbox{ let } \phi(v)=B(v)\cdot \sigma,
% \forall v \in \r:
% \\ 
% \{ 
%\phi'(v)=\cdots=\phi^{(\gamma-1)}(v)=0 
%\}
% \Longrightarrow \phi^{(\gamma)}(v)\neq 0.
%  \end{array} \qquad (\gamma ND)
%$$
 \begin{eqnarray*}
 \forall (v,\sigma) \in  \r^M  \times S^{N}, 
 &\qquad 
  \displaystyle{ \sum_{k=0}^{\gamma -1} \left| D^k b(v)\cdot \sigma \right| }
  > 0.
&\qquad (\gamma ND)
 \end{eqnarray*}

\noindent Before  proving the Theorem  \ref{FODELPNS}
   we give  some useful results about oscillatory  integrals
  following  Stein's book \cite{Ste}.
\begin{Prop}[\cite{Ste}]\label{IntExp}
%-......................................
%$\mbox{}$ //
   Suppose  $\phi \in C^{k+1}(\r,\r)$
  so that, for some $k \geq 1$, 
 \begin{eqnarray}  \label{HypPhi2}
  \dis \frac{d^k\phi}{d v^k}(v)  \geq 1,   &\qquad&  \forall v \in ]\alpha,\beta[. 
 \end{eqnarray}
Then 
$$\dis 
 \left| \int_\alpha^\beta e^{i \lambda \phi(v)}  dv \right| 
 \leq c_k \cdot \frac{1}{|\lambda|^{1/k}}
$$ 
holds when 
  \begin{enumerate} 
   \item $k \geq 2$ or 
  \item $ k=1$ and $\phi'$ is monotonous. 
\end{enumerate}   
 Furthermore, the bound  $ c_k$ is independent of $\lambda$ and $\phi$.
 \end{Prop}
 This Proposition can be found in 
 \cite{Ste} p 332.  
 Elias  M. Stein obtains $ c_k \leq 5 \cdot 2^{k-1}-2$
 in his proof. Notice that $c_k$ is independent of the length of the interval
   $]\alpha,\beta[$. 
  For $|\lambda| < 1$, the  bound for the oscillatory integral blows up. 
  Indeed, for $k=1$,  we can relax the monotonous assumption on $\phi$  
 by the following  bounds 
\begin{eqnarray*}  \label{HypPhilowerbound}
  \dis |\phi'(v)|  \geq \delta > 0, \quad  \forall v \in ]\alpha,\beta[,
  & \qquad \ds
    \widetilde{c}_1= 2 + \delta^{-1}\int_\alpha^\beta |\phi"(v)|dv, 
 \end{eqnarray*}
 Indeed, integrating by parts and using the inequality
 \\ 
  $\min(a,\beta b) \leq \min(1,\beta)\max(a,b)$ 
 for all non negative $a,b,\beta$, we get
 $$
 \left| \int_\alpha^\beta e^{i \lambda \phi(v)}  dv \right| 
 \leq \max(|\beta-\alpha|,\widetilde{c}_1) \cdot \max(1,\frac{1}{\delta})  \cdot \min(1,\frac{1}{|\lambda|}).
$$ 
\\
Furthermore, the bound given in Proposition \ref{IntExp} blows up for small 
 $\lambda$, so we replace it  by the length of the interval and get 
 the following Corollary.

 \begin{Cor}\label{CorUIntExp}
%-......................................
%$\mbox{}$ //
  Let  $\delta > 0$. Suppose  $\phi \in C^{k+1}(\r,\r)$
  so that, for some $k \geq 1$, 
 \begin{eqnarray}  \label{HypPhi}
  \dis \left|\frac{d^k\phi}{d v^k}(v) \right|  \geq \delta, 
  &\qquad&  \forall v \in ]\alpha,\beta[. 
 \end{eqnarray}
Then 
$\dis 
 \left| \int_\alpha^\beta e^{i \lambda \phi(v)}  dv \right| 
 \leq \max(|\beta-\alpha|,\widetilde{c}_k) \cdot \max(1,\frac{1}{\delta^{1/k}})
  \min(1,\frac{1}{|\lambda|^{1/k}}),
$
\\ 
where $ \widetilde{c}_k$ is independent of
  $\lambda$, $\phi$ and $]\alpha,\beta[$ for $k \geq 2$ 
 \\ 
and 
 $ \dis  \widetilde{c}_1= 2 + \delta^{-1}\int_\alpha^\beta |\phi"(v)|dv$.
 \end{Cor}
Notice that, for $k \geq 2$, 
 $ \widetilde{c}_k =c_k$ is given in Proposition \ref{IntExp}.
%  and the constant  $\min(|\beta-\alpha|,\widetilde{c}_k)$ is enough.
\\
Following Stein's book (Corollary p 334), 
we obtain the  following Proposition.
 \begin{Prop}[\cite{Ste}]\label{LPNSU0}
%-......................................
%$\mbox{}$ //
   Let $\psi \in W^{1,1}(]\alpha,\beta[)$,
$\phi \in C^{k+1}(\r,\r)$ such
 that, for some  $\delta > 0$  and $ k \geq 1$, 
\bqre
  \dis \left | \frac{d^k\phi}{d v^k}(v)\right| & \geq \delta, \quad \forall v \in 
]\alpha,\beta[.
  \eqre
Then
$$\dis 
 \left| \int_\alpha^\beta \psi(v) e^{i \lambda \phi(v)}  dv \right| 
 \leq 
 \frac{\max(|\beta-\alpha|,\widetilde{c}_k)}{ \min(1,\delta^{1/k})\max(1,|\lambda|^{1/k}) )}
  \left( 
          \|\psi\|_{L^\infty(]\alpha,\beta[)}  +  \|\psi'\|_{L^1(]\alpha,\beta[)}
  \right),
$$ 
where $ \widetilde{c}_k$ is independent of $\lambda$, $\phi$, $\psi$
 and $]\alpha,\beta[$ for $k \geq 2$,
\\
and $\ds  \widetilde{c}_1= 2 + \delta^{-1}\int_\alpha^\beta |\phi"(v)|dv$.
 \end{Prop}
\pr
 This is classically proved in writing the integral
 $ \dis  \int_\alpha^\beta \psi(v) e^{i \lambda \phi(v)}  dv $
as   $ \dis  \int_\alpha^\beta \psi(v)I'(v)  dv $, 
 with 
  $ \dis I(v) =  \int_\alpha^v e^{i \lambda \phi(u)}  du $,
 integrating by parts and using the uniform  estimate 
   for $|I(v)|$ 
 % \leq \max(|v-a|,\widetilde{c}_k) \cdot \delta^{-1/k}(1+|\lambda|)^{-1/k}$
 from previous Corollary.
$\CQFD$

\bigskip

%\begin{Rk}: \label{RLemmaStein}
% The Proposition \ref{LPNSU0} can be found in a Elias,  M. Stein's book, 
% \cite{Ste} p.342, with $\psi$ smooth.
% Reading the proof therin, we can see that 
%$\psi \in W^{1,1}(\r)$  with compact support 
% is enough to give the same result.
 % Furthermore, if the lower bound of $|\phi^{k}|$ is replaced by  a positive 
% constant $\delta > 0$,  this only changes $ \widetilde{c}_k(\phi)$ by 
% $ \widetilde{c}_k(\phi)/ \delta^{1+1/k}$.  In  a similar way, 
% if the length of the support of $\psi$  is less than $l>0$, it suffices
% to  change 
% $\widetilde{c}_k(\phi)$ by  $l^{k+1} \widetilde{c}_k(\phi)$. 
%   These last remarks are useful in the proof of  the Theorem  \ref{FODELPNS}.
%\end{Rk}
%
 Now we generalize Proposition \ref{LPNSU0} in the case with parameters 
and a  $(\gamma ND)$ like assumption.
\begin{Prop}\label{LPNSU1}
%-......................................
%$\mbox{}$ //
   Suppose  $P$ is a compact set of parameter $p$, $A > 0$,
 $\psi(u;p)$  belongs to $L^\infty_p(P,W^{1,1}_u(]-A,A[))$
 and $\phi(u;p) \in C^{\gamma +1}(\r_u\times P_p, \r)$, 
 such that, for all $(u,p)$ in 
$K= [-A,A] \times  P$,
 \bqr 
 \label{cond}
\dis 
  \sum_{k=1}^{\gamma} 
 \left| \frac{\partial^k \phi}{\partial u^k} \right| (u;p)& >& 0.
\eqr
Then, for any $]\alpha,\beta[ \subset ]-A,A[$,  
\bqre 
 &  \dis 
  \left| \int_\alpha^\beta \psi(u;p) e^{i \lambda \phi(u;p)}  du \right| 
 \\
 \leq 
 & d_\gamma \cdot \min \left(1,\ds \frac{1}{|\lambda|^{1/\gamma}} \right) \cdot
   \left( \|\psi\|_{L^{\infty}(K)} + 
 \left \|\frac{\partial \psi}{\partial u} \right\|_{L^\infty(P,L^{1}(]-A,A[))} \right),
\eqre
where 
 constant $ d_\gamma$ is independent of $\lambda$ and 
only depends on  $A$, 
$ \dis \sup_{K}\left | \frac{\partial^2 \phi}{\partial u^2}\right | $,
$ \dis \inf_{K}  \frac{1}{\gamma} \sum_{k=1}^{\gamma} 
 \left| \frac{\partial^k \phi}{\partial u^k} \right| $.
  % and $\phi$ restricted to $K$.  
 \end{Prop}
%*..........................................
\pr 
 Since $K$ is a compact set, we can choose $0 < \delta\leq 1$ 
 such that, everywhere on $K$: 
$$ 
 0 < \delta < \frac{1}{\gamma} \sum_{k=1}^{\gamma} 
 \left| \frac{\partial^k \phi}{\partial u^k} \right| (u;p).
$$
Let us define
 the open set 
$  Z_k =\{ (u;p), \, |\partial_u^k\phi (u;p)|> \delta \}$,
 for $k=1,\cdots,\gamma$. 
Necessarily $ \displaystyle{ K \subset \bigcup_{k=1}^{\gamma} Z_k}$, 
and then there exists a partition of unity such that
 $ \displaystyle{ \sum_{k=1}^{\gamma  } \rho_k \equiv 1}$ on $K$ and such that 
  the support of $\rho_k$ is included in $Z_k$.
 Let us define $\psi_k =~\rho_k \psi$ and  $I= I_1+\cdots + I_\gamma$
 where 
$  \displaystyle{ I_k(p)=
   \int_a^b \psi_k(u;p) e^{i \lambda \phi(u;p)}  du  } $.
 We   apply Proposition \ref{LPNSU0} 
 on each $I_k$ where the exponent `` $'$ '' denotes $\partial_u$:
 
 $$  |I_k| \leq  \frac{\max(2A,\widetilde{c}_k)}{\delta^{1/k} 
   \max( 1,|\lambda|^{1/k})} 
            \sup_P\left (\| \psi_k(.,p) \|_{L^\infty(]-A,A[) } +
                   \| \psi_k'(.,p) \|_{L^1(]-A,A[) }  \right).  $$ 
Since for any fixed $p$ and $J=]-A,A[$, we have
    \begin{eqnarray*}
&& \ds      \left (\| \psi_k(.,p) \|_{L^\infty(J) } +
                   \| \psi_k'(.,p) \|_{L^1(J) }  \right)
 \\ &\leq  & \ds
  \left (\| \rho_k \|_{L^\infty(J) } +
                   \| \rho_k' \|_{L^1(J) }  \right)
\left (\| \psi(.,p) \|_{L^\infty(J) } +
                   \| \psi'(.,p) \|_{L^1(J) }  \right), 
\end{eqnarray*}  
it is enough to take
 $$ \dis 
 d_\gamma = 
            \sum_k  \frac{\max(2A,\widetilde{c}_k)}{\delta^{1/k}}
    % \max( 1,|\lambda|^{1/k}) 
            \left (\| \rho_k \|_{L^\infty(K) } +
                   \| \partial_u \rho_k \|_{L^\infty(P,L^1(J)) }  \right)
 $$
 to  conclude the proof.
$\CQFD$
%*....................................................

\bigskip

\noindent We are now able to prove  the  second Theorem.\\

\noindent {\bf Proof of  Theorem  \ref{FODELPNS}.}\\
   The proof is splitted in three steps.
 First, we choose a suitable variable associated to $D$.
 Secondly, we use Fourier transform with respect to $X$ and  solve 
 a linear  ordinary differential equation with respect to $v_1$. 
 Third, we  obtain Sobolev estimates for $\rho_\psi$
  with Proposition \ref{LPNSU1}.  
\\

\underline{Step 1, change of coordinates: }
  With a suitable choice of orthonormal coordinates, we assume,
 without loss of generality that  
   $$
   D= F \cdot \nabla_v  = |F| \frac{\partial}{\partial v_1}
   $$ 
 where $|F|$ is the euclidean norm of vector $F$ and 
 $ v = (v_1,v_2,\cdots,v_M) \equiv (v_1;w)$. 
  Notice that the jacobian for an orthonormal change of variables is one,
  thus the estimates on $\rho_\psi$   are invariant through such choice for $v_1$.
  With such notations, equation (\ref{eqb})
  becomes
  \bqr  \label{eqv1}
%$$$$$$$$$$$$$$$$$$$$$$$$$$$$$$$$$$$$$$$$$$$$$$$$$$$$$$
 b(v) \cdot \nabla_X f +|F|\frac{\partial f}{\partial v_1} =g.
\eqr

\underline{Step 2, linear o.d.e.: }
   Denoting by $ {\cal F}(f)$ the Fourier transform of $f$ with respect to $X$,
   and by $Y$ the dual variable of $X$,
   equation (\ref{eqv1}) becomes
  \bqr 
   \label{ode1}
   |F| \frac{\partial }{\partial v_1}{\cal F}(f)  +i (b(v)\cdot Y) {\cal F}(f) 
  &=& {\cal F}(g).
  \eqr
 For almost all fixed $Y$, 
we solve an ordinary differential equation with respect 
 to $v_1$. For this purpose, 
  we chose the initial $v_1$, namely $v_1^0 \in ]0,1[$,  such that 
 \bqr \label{Fubini0} 
    \begin{array}{cl}
    & \dis 
  \int_{\r^{N+1}_Y \times \r^{M-1}_w} |{\cal F}(f)|^2(Y;v_1^0;w) dY dw  \\
    \leq &
  \dis     \int_{\r_{v_1}} 
   \int_{\r^{N+1}_Y \times \r^{M-1}_w} |{\cal F}(f)|^2(Y;v_1;w) dYdwdv_1. 
 \end{array}
 \eqr
 Existence of such $v_1^0$ is a consequence of Fubini's Theorem. 
 \\ Indeed, let   $ \displaystyle{
                h(v_1)= 
  \int_{\r^{N+1}} \int_{\r^{M-1}_w} |{\cal F}(f)|^2(Y;v_1;w) dYdw    \geq 0
                    }$.
\\
 Function $h$ is defined almost everywhere,
  belongs to $L^1(\r_{v_1})$ and satisfies 
   $\|h\|_{L^1(\r_{v_1})} 
      =\|f\|^2_{L^2_{X,v}}  $.
 Since $h$  function cannot be everywhere greater than its mean value
 on $]0,1[$,
 there exists $v_1^0 \in ]0,1[$
 such that
  $ \displaystyle{
        h(v_1^0) \leq 
       \int_{0}^1    h(v_1) d v_1       
            },$
which confirms (\ref{Fubini0}).
%the  existence of such $v_1^0$.
 % \\
 %\medskip
%

\noindent We finally write an explicit formula for $ {\cal F}(f)$
 with  $B(v)$ being a primitive with respect to $v_1$
of $-b/|F|$: 
\bqre
 B(v) & = & B(v_1;w) = -\int_{v_1^0}^{v_1} \frac{ b(u;w)}{|F|} du 
\\
{\cal F}(f)(Y,v_1;w) &=& {\cal F}(f)(Y,v_1^0;w) e^{ i  B(v)\cdot Y }
   \\ && 
  + \frac{1}{|F|} \int_{v_1^0}^{v_1}  {\cal F}(g) (Y,u;w)
  e^{i  (B(v_1;w)-B(u;w))\cdot Y} du.
\eqre

\underline{Step 3, $H^{1/\gamma}$ estimates with oscillatory integrals}:
We decompose 
$\dis \rho_\psi(t,x)=\int_{\r^M} f(t,x,v) \psi(v)\,dv$ 
   in two parts from the explicit expression of ${\cal F}(f)$ in step 2:
${\cal F} (\rho_\psi)=\widehat{\rho}_f + \widehat{\rho}_g$.
The first term is
\bqre
   \widehat{\rho}_f(Y) & =& \int_{\r^{M-1}_w}{\cal F}(f)(Y,v_1^0;w) 
             \int_{\r_u} \psi(u;w) e^{ i  B(u;w)\cdot Y }du dw.
 \eqre
In this integral, there is an oscillatory integral which is parametrized by 
$w$ and 
$Y=\lambda \sigma $ with $\lambda =|Y|$ and $\sigma \in S^N$; it is 
 \bqr  \label{oscI} 
      Osc(Y,w) & = & 
    \int_{\r_u} \psi(u;w) e^{ i \lambda  B(u;w)\cdot \sigma }du.
 \eqr
 To use  the  Proposition \ref{LPNSU1},  we set 
  $p=(\sigma,w)$ which belongs  to the compact set 
$ P= S^N\times [-A,A]^{M-1}$ 
  with $A > 1 > v_1^0 > 0$ such that $ supp \, \psi \subset  [-A,A]^M$.
Condition (\ref{cond}) of  Proposition \ref{LPNSU1} for 
 oscillatory integral (\ref{oscI}) is
 \bqre 
   \dis  \sum_{k=1}^\gamma \left| \frac{\partial^k  B(u;w)}{\partial u^k}
                           \cdot \sigma \right|
  &> 0
 \eqre
which is exactly the  $(\gamma ND)$ assumption for $b(.)$.
Thanks to the $(\gamma ND)$ assumption and 
 Proposition \ref{LPNSU1}, there exists a constant L 
 such that for all $( Y,w) \in \r^d\times[-A,A]^{M-1},$
    and for all $\alpha, \beta$ such that   $-A < \alpha < \beta < A$,
  we have
 \bqr \label{L} 
     \max(1,  |Y|^{1/\gamma})|
   \left|
  \int_\alpha^\beta  \psi(u;w) e^{ i \lambda  B(u;w)\cdot \sigma }du
        \right|  
  & \leq & L.
\eqr
Using  constant $L$ and the compact support of $\psi$ we have
 \bqre 
      \max(1,  |Y|^{1/\gamma})|\widehat{\rho}_f(Y)|  & \leq & 
    L  \int_{[-A,A]^{M-1}} 
    |{\cal F}(f)(Y,v_1^0;w)|dw .
 \eqre
By  Cauchy-Schwarz inequality, we get
\begin{eqnarray*}
 \max(1,  |Y|^{2/\gamma})|\widehat{\rho}_f(Y)|^2 &\leq  &
     (2A)^{M-1} L^2 \int_{[-A,A]^{M-1} } 
    |{\cal F}(f)(Y,v_1^0;w)|^2dw
.  
\end{eqnarray*}
Finally, since $v_1^0$ satisfies (\ref{Fubini0}), we obtain
\begin{eqnarray*}
 \int_{\r^{N+1}}
 \max(1,  |Y|^{2/\gamma})|\widehat{\rho}_f(Y)|^2 dY  &\leq  &
     (2A)^{M-1} L^2 \int_{\r^{N+1}\times\r^{M} } 
    |{\cal F}(f)(Y,v)|^2 dv dY, 
\end{eqnarray*}
which gives 
             $\widehat{\rho}_f \in H^{1/\gamma}$.
\bigskip 

\noindent  The second term  $\widehat{\rho}_g$ is bounded in the same way.
  More precisely, we set
 \bqre
  \widehat{\rho}_g(Y)  &= &   \int_{\r^{M-1} }  H(Y,w) dw  % \mbox{ with: }
  \eqre
with 
 \bqre
  H(Y,w) & = &
  \dis  \frac{1}{|F|} \int_{-A}^{A} \int_{v_1^0}^{v_1} 
       {\cal F}(g)(Y,u;w) e^{ i  (B(v_1;w)-B(u;w))\cdot Y }du  dv_1 . 
   \eqre
Using   Fubini's Theorem 
 and   notation  
 \bqre 
   \Psi(Y,u;w) &=&  \psi(u;w) e^{ i  B(u;w)\cdot Y},
 \eqre
  we have  another expression of
 $H(Y,w)$:
\bqre 
 H(Y,w)  &= & 
     \dis \frac{1}{|F|}
    \int_{v_1^0}^{A} 
  {\cal F}(g)(Y,u;w) e^{- i  B(u;w)\cdot Y}
      \left( \int_{u}^{A} 
    \Psi(Y,v_1;w)  dv_1\right) du 
\\  & &  
  + \dis    \frac{1}{|F|}  \int_{-A}^{v_1^0} 
  {\cal F}(g)(Y,u;w) e^{- i  B(u;w)\cdot Y} 
   \left( \int_{-A}^{u} \Psi(Y,v_1;w)dv_1\right)du, 
\eqre 
where there are two oscillatory integrals 
 $\dis \int_{u}^{A} \Psi(Y,v_1;w) dv_1 $ and  $\dis \int_{-A}^{u} \Psi(Y,v_1;w) dv_1 $
which are uniformly bounded thanks to inequality (\ref{L}).
Then we have
 \bqre
\max(1,|Y|^{1/\gamma})|H(Y;w)| 
 & \leq  &
\dis     \frac{L}{|F|} \int_{-A}^{A} \left|{\cal F}(g)(Y,u;w) \right|  du.
\eqre
With  Cauchy-Schwarz inequality, we obtain
\bqre
 \max(1,|Y|^{2/\gamma})|H(Y;w)|^2 
    & \leq  &
\dis     \frac{2A L^2}{|F|^2}
 \int_{-A}^{A} \left|{\cal F}(g)(Y,u;w) \right|^2  du
\eqre
and finally
\bqre
\max(1,|Y|^{2/\gamma}) | \widehat{\rho}_g(Y) |^2 
 & \leq &
 (2A)^M    \frac{L^2}{|F|^2} \int_{\r^{M}}
   |  {\cal F}(g)(Y,v)|^2 dv.
 \eqre
Then $ \rho_g \in H^{1/\gamma}$, thus finally $\rho_\psi$ is also in this space,
which concludes the proof of the Theorem.
\CQFD

%%%%%%%%%%%%%%%%%%%%%%%%%%%%%%%%%%%%%%%%%%%%%%%%%%%%%%%%%%%%%%%%%%%%%%%%%%%%%%%
%%%%%%%%%%%%%%%%%%%%%%%%%%%%%%%%%%%%%%%%%%%%%%%%%%%%%%%%%%%%%%%%%%%%%%%%%%%%%%
\Section{About  non degeneracy conditions}\label{scomparisons}
%%%%%%%%%%%%%%%%%%%%%%%%%%%%%%%%%%%%%%%%%%%%%%%%%%%%%%%%%%%%%%%%%%%%%%%%%%%%%%
%%%%%%%%%%%%%%%%%%%%%%%%%%%%%%%%%%%%%%%%%%%%%%%%%%%%%%%%%%%%%%%%%%%%%%%%%%%%%%%

 Theorem \ref{lemmemoyenne} and Theorem \ref{FODELPNS} assume two different 
  non degeneracy conditions on vector field  $a(v) \in \r^N$, 
$v \in \mbox{ supp } \psi \subset \r^M$. Those conditions involve two parameters,
 namely $\alpha=\alpha_{a(.)} \in ]0,1]$ in (\ref{alpha}) 
 and $\gamma = \gamma_{a(.),F} \in \N^*$ in  (\ref{gamma}),
 directly linked to the smoothing effect for the averaging in 
 $H^{\alpha/2}_{loc}$ or  $H^{1/\gamma}_{}$.
  In this section, we give some optimal upper bounds for
 $\alpha$ and $1/\gamma$ to compare 
both results
 obtained by different  ways.
     Indeed, for $M=1$ and $N\geq 2$,
 Theorem  \ref{FODELPNS} gives a better smoothing 
  effect than Theorem \ref{lemmemoyenne}.
  Conversely, when $N=M$, Theorem \ref{lemmemoyenne} 
 is stronger than Theorem \ref{FODELPNS}. 
In this part, we study these various properties and in particular, we prove Theorem 3.
\\
More precisely, let $A$ be positive, 
 we obtain the optimal $\alpha$ and $\gamma$, namely
\bqre 
 \alpha_{opt}(N,M)  & = & \sup_{a(.) \in C^\infty([-A,A]_v^M,\r_x^N)} \alpha,
 \\ % & \quad &
 \gamma_{opt}(N,M) & =& 
 \min_{a(.) \in C^\infty(\r_v^M,\r_x^N),\, F \in \r^N \setminus\{0\}} \gamma.
\eqre
%$  \alpha_{opt}(N,M) $
% is the supremum of $\alpha$,
%$ \gamma_{opt}(N,M) $
% is the minimum of $\gamma$ for all $a(.) \in C^\infty(\r_v^M,\r_x^N)$.
%
 We start by obtaining the easiest estimate which is a lower bound for $\gamma$.
\begin{Prop}
 \label{mingamma}
 For all $N, M$, we have $\gamma \geq \gamma_{opt}(N,M)= N+1.$
\end{Prop}
{\bf Proof.}
 We use notations from Section \ref{sF}.
Following this section,
the $ (\gamma ND)$ condition can be rewritten and 
 means that we cannot find  
 $\sigma \in S^{N}$ such that 
 $\sigma \perp b(v)$, $\sigma \perp D b(v)$, \ldots, 
$\sigma \perp D^{\gamma -1}b(v)$. 
    There are $\gamma$  conditions to satisfy. 
  Since $b(v)$ belongs to $\r^{N+1}$,  we necessarily  have 
 %\bqr
  $
\gamma \geq  N+1.
$
% \eqr
 Indeed $N+1$ is the minimal possible value for $\gamma$.
%\\
 For instance,  if
  $D =  \displaystyle{\frac{\partial  }{\partial v_1}}$, 
    $b(v)=(1,v_1,v_1^2,\cdots,v_1^{N})$,
  with $v=(v_1,v_2,\cdots,v_M)$, 
 we have $\gamma_{opt}=N+1$.
%\\
% A discussion and  links 
% with the classical Lions-Perthame-Tadmor non-degeneracy condition  
% is done in \cite{Ju}. 
$\CQFD$
\medskip \\
The optimal $\alpha$ is more difficult to get and it is obtained 
in the following   subsections, see also \cite{Ju}. 
 The evaluation of   exponent $\alpha$  also implies 
 new asymptotic expansions involving  piecewise smooth functions
  in  \cite{JR1}.

%%%%%%%%%%%%%%%%%%%%%%%%%%%%%%%%%%%%%%%%%%%%%%%%%%%%%%%%%%%%%%%%%%%%%%%
\Subsection{$M=1$, one dimensional velocity} \label{ssM=1}    
%%%%%%%%%%%%%%%%%%%%%%%%%%%%%%%%%%%%%%%%%%%%%%%%%%%%%%%%%%%%%%%%%%%%%%%
\begin{Prop}
 \label{maxalpha}
 For $M=1$, we have $\alpha \leq \alpha_{opt}(N,1)= \dis \frac{1}{N}.$
\end{Prop}
%{\bf Proof :}
%
%$\CQFD$
%\medskip \\
 To obtain this optimal $\alpha$ for $M=1$, we need some other 
notations and the following results. The proof of Proposition \ref{maxalpha}
 is achived at the end of this subsection \ref{ssM=1}.
\medskip\\
Let $\ph \in C^\infty([a,b],\r)$ and $v \in [a,b]$, the multiplicity of $\ph$ 
 on $v$ is defined  by 
\bqre 
 m_\ph[v] &= & \inf \{ k \in \N,\, 
       \ph^{(k)}(v) \neq 0 \}   
 \qquad   
 \in \overline{\N}= \N \cup \{+\infty\}.
\eqre
It means that if $k=m_\ph$ then $\ph^{(k)}(v) \neq 0$ and 
  $  \ph^{(j)}(v)= 0 \mbox{ for }  j=0,1,\cdots,k-1        
$.
For instance $m_\ph[v]=0$ means $\ph(v) \neq 0$;
   $m_\ph[v]=1$ means $\ph(v)=0$, $\ph'(v) \neq 0$
  and $m_\ph[v]=+\infty$ means $\ph^{(j)}(v)=0$ for all $j \in \N$.
\\
 Set the multiplicity of $\ph$ on $[a,b]$ by
\bqre
 m_\ph  & = & \sup_{v \in [a,b]} m_\ph[v]  \quad  \in \overline{\N}.
\eqre
Notice that the case where $\ph$  only belongs to $C^k$, $m_\ph$ is well defined only
 if $m_\ph[v]\leq k$ for all $v \in [a,b]$.

\begin{Lemma} 
 \label{mphi}
 Let $\ph \in C^k([a,b],\r)$ with $a < b$, and 
 \bqre 
    Z(\ph,\eps)  & = & \{ v \in [a,b],\; |\ph(v)|\leq \eps\}.
 \eqre
If $m_\ph$ is well defined ($m_\ph \leq k$)
then there exists $C> 0 $ such that, for all $\eps > 0$,
  \bqr  \label{mphineq}
    \meas(Z(\ph,\eps))\leq  C \eps^\alpha 
  & \mbox{  with  } & 
 \dis \alpha =\frac{1}{m_\ph}.
 \eqr
Furthermore, if $m_\ph$ is positive, for all $\beta > \alpha$, we have
  $ \dis 
      \lim_{\eps \rightarrow 0}\frac{ \meas(Z(\ph,\eps))}{\eps^\beta}
   = + \infty $ (Optimality).
\end{Lemma}
{\bf Proof.} 
The case $m_\ph=0$ is clear enough since there is no zero in this situation.
Quantity $m_\ph$ is positive simply means that the set $Z(\ph,0)$ of roots of
 $\ph$ is not empty. Since any root of $\ph$ has a finite multiplicity, 
  the compact set  $Z(\ph,0)$ is discrete and then finite: 
 $Z(\ph,0)=\{z_1,\cdots,z_\nu\}$.  For each $z_i$ and $h>0$, 
 let $V_i(h)$ be $]z_i-h,z_i+h[ \cap [a,b]$. For any $0<h<|b-a|$, we have 
$$h \leq \mbox{ meas} (V_i(h)) \leq 2h.$$
%\\
For any root $z_i$,  there exists $h_i \in ]0,|b-a|[$,
  $A_i > 0$ 
  and $\delta_i>0$ such that
\bqr \label{eqpolyk}
 \dis   \delta_i |h|^{k_i} &\leq  |\ph(z_i+h)|
 \leq &  A_i |h|^{k_i} \quad  \mbox{ for all }  h \in V_i(h_i),
\eqr
with $k_i=m_\ph[z_i]$. This is a direct consequence of
Taylor-Lagrange formula.
 Let $\dis V$ be $ \dis \bigcup_{i}V_i(h_i)$ and 
 $\eps_0 = \min\left(1, \mind_{v\in [a,b] \setminus V}|\ph(v)| \right)$.
By the continuity of $\ph$ on the compact set $[a,b] \setminus V $, 
$\eps_0$ is positive.
Then for all $0<\eps < \eps_0$, we have 
 $Z(\ph,\eps) \subset V$.
If $\eps \geq  |\ph(z_i+h)|$ for $|h|<h_i$, 
then from (\ref{eqpolyk}), we have  $(\eps/\delta_i)^{1/k_i} \geq |h|$.
This last inequality 
 implies  for $0<\eps < \eps_0 \leq 1$  that 
 $Z(\ph,\eps)$ is a subset of $\dis \bigcup_i V_i((\eps/\delta_i)^{1/k_i})$ 
and then 
$$  
    \mbox{meas} (Z(\ph,\eps)) 
    \leq \dis 2 \sum_{i=1}^\nu (\eps/\delta_i)^{1/k_i} 
  \leq \left( 2\sum_{i=1}^\nu \delta_i^{-1/k_i}  \right)  \eps^{1/m_\ph}.
$$
It gives inequality (\ref{mphineq}).
To obtain the optimality of $\alpha$, let $z_j$ be a root of $\ph$
 with maximal multiplicity i.e.  $m_\ph[z_j]=m_\ph =k$.
Again  from (\ref{eqpolyk}),  
$ V_j((\eps/A_j)^{1/k})$ is a subset of $Z(\ph,\eps)$ for  all 
$ \eps \in]0,\eps_0[$.
Then we have  
$
(\eps/A_j)^{1/k}\leq  \mbox{meas} (Z(\ph,\eps)) 
$, which is enough to get the optimality of $\alpha=1/k$
and  concludes the proof.
$\CQFD$
\medskip \\
 An upper bound of $\alpha_{opt}(N,1)$ is  a consequence of previous 
  Lemma.
\begin{Lemma}\label{LsinfN}
  For all $N$, we have $\alpha_{opt}(N,1)\leq 1/N$.
\end{Lemma}
\pr 
 For any  $a(.) \in C^\infty(\r_v,\r_x^N)$ and  $A>0$,
we set
\bqre 
 \ph(v;u,\sigma) & = & a(v) \cdot \sigma -u= b(v)\cdot (-u,\sigma),
 \eqre 
 defined for $v \in [-A,A]$,
 with $ u \in \r$, 
  $ (-u,\sigma) \in S^{N}$, 
  $b(v)=(1,a(v))\in \r^{N+1}$
 and 
 $m = \dis \sup_{ (-u, \sigma) \in S^{N}} m_{\ph(.;u,\sigma)}$. 
 \\  %
  Let $v$ be fixed, we choose $(-u,\sigma)$ such that $m_\ph[v]\geq N$
 in order to  obtain a lower bound for $m$.\\
 Since $rank\{b(v),b'(v),\cdots,b^{(N-1)}(v)\} \leq N$, there exists 
   $(-u,\sigma)$  such that $ u^2+|\sigma|^2=1$ and 
$ (-u,\sigma)\perp \{b(v),b'(v),\cdots,b^{(N-1)}(v)\}$.
 Then with such $u$ and $\sigma$,  $m_{\ph(.;u,\sigma)}[v]\geq N$
  which implies $m \geq N$ and consequently,
 from the optimality  obtained in Lemma \ref{mphi}, we get 
 $\alpha \leq \alpha_{opt}(N,1) \leq  \dis \frac{1}{N}.$
$\CQFD$
\medskip \\
When function $v \rightarrow \ph(v;p)$ depends on a parameter $p$,
some results  are obtained in the two following Lemma to bound quantity $C$ of 
 Lemma \ref{mphi} independently of $p$ parameter.

\begin{Lemma}\label{Ls2}
Let $k \geq 1$, $I$ an interval of $\r$,
    $ \phi \in C^k(I,\r)$ and $\delta > 0$.
\\
If $|\phi^{(k)}(v)| \geq \delta > 0$ for all $x \in I$
then there exists a constant $\overline{c}_k$ independent of $\phi, I, \delta$
such that 
\bqre  \mbox{meas}( Z(\phi,\eps)) \leq \overline{c}_k (\eps/\delta)^{1/k}, 
\quad \mbox{ where }   
        Z(\phi,\eps)=\{v \in I,\;  |\phi(v)|\leq \eps\}.
\eqre
\end{Lemma}
{\bf Proof.} Since the result is independent of  interval $I$ 
and of  $\phi^{(k-1)}(0)$ sign,
 let us suppose that $I=\r$ with $|\phi^{(k)}(v)| \geq \delta > 0$ on $\r$,
  and $\phi^{(k-1)}(0)\leq 0$.\medskip
\\
 We first treat the case $k=1$.
If $\phi'(v)$ stays positive, we have $\phi(0)+\delta v \leq \phi(v)$ for $0 \leq v$ and
 since $\phi(0)\leq 0$, there exists a unique  $c \geq 0$
 such that $\phi(c)=0$. 
In the other case, $\phi'(v)$ stays negative, 
and we find a unique  $c \leq 0$
 such that $\phi(c)=0$. 
%Changing $v$ by $v+c$, let us suppose that $c=0$.
% With this translation, the main assumption on $\phi$ 
%and  the measure of $Z(\phi,\eps)$ is invariant.
 Then $|\phi(v)|\geq  \delta |v-c|$ for all $v$, and $|\phi(v)| \leq \eps$
 implies $|v-c|\leq \eps/\delta$ i.e. 
  $Z(\phi,\eps) \subset [c-\eps/\delta,c+\eps/\delta]$.
So the lemma is proved for $k=1$ with $\overline{c}_1=2$.  
\medskip\\
We now prove the Lemma by induction on $k$. 
Let us suppose that the case $k$ is known.
As for $k=1$, there exists 
 a unique $c$ such that $\phi^{(k)}(c)=0$. 
 Thus for all $v$ we have $|\phi^{(k)}(v)| \geq \delta|v-c|$.
Let $\eta> 0 $ and set 
    $W=Z(\phi,\eps) \cap [c-\eta,c+\eta]$, 
    $U=Z(\phi,\eps) \cap (]-\infty,c -\eta[ \cup ]c+\eta,+\infty[)$.
We have
$\mbox{meas}(W) \leq 2 \eta$ and  by our inductive hypothesis, since
$|\phi^{(k)}(v)| \geq \delta|v-c| \geq \delta \eta$ on $U$,
$\mbox{meas}(U) \leq \overline{c}_k (\eps/(\delta \eta))^{1/k}$.
Now the relation $Z(\phi,\eps)=W \cup U$ gives   
$ \mbox{meas}( Z(\phi,\eps)) \leq 
  \dis \inf_{\eta> 0} \left(2\eta + \overline{c}_k (\eps/(\delta \eta))^{1/k}\right)$
 which
 implies by a simple computation of the minimum that
 $
   \mbox{meas}( Z(\phi,\eps)) \leq 
  \dis  \overline{c}_{k+1}(\eps/\delta)^{1/(k+1)}$,
 where $  \overline{c}_{k+1}=2^{1/(k+1)}(k+1)k^{1/(k+1) -1 }\overline{c}_k^{1-1/(k+1)}   $  
which concludes the proof.
$\CQFD$

\begin{Lemma}\label{Lscompact}
 Let $P$ be a compact set of parameters, 
$k$ a positive integer, 
$A > 0$, $V=[-A,A]$,  $K= V \times  P$,
 $\phi(v;p) \in C^0(P,C^{k}(V, \r))$, 
 such that, for all $(v,p)$ in the compact $K$, we have 
\bqre 
 \label{cond2}
\dis 
  \sum_{j=1}^{k} 
 \left| \frac{\partial^j \phi}{\partial v^j} \right| (v;p)& >& 0.
 \eqre
 Let $ Z(\phi(.;p),\eps)=\{v \in V,\;| \phi(v;p)| \leq \eps\}$, 
then  there exists  a constant $C$ such that
\bqre   
 \dis  \sup _{p\in P} \mbox{meas}(Z(\phi(.;p),\eps)) &\leq &C \eps^{1/k}.
\eqre
\end{Lemma}
{\bf Proof.}
 Since $K$ is a compact set, we can choose $0<\delta\leq 1$ 
 such that, everywhere on $K$, we have
$ \dis  
 0 < 2\delta < \frac{1}{k} \sum_{i=1}^{k} 
 \left| \frac{\partial^i \phi}{\partial v^i} \right| (v;p).
$
\\
For each $(v;p) \in K$, there exists an integer $i\in \{1,\cdots,k\}$,
  a number $r>0$ and an open set $O_p$ with 
 $p \in O_p \subset P$ such that 
$|\partial^i_v \phi|>\delta$ on $U(v,p)=]v-r,v+r[\times O_p$.
Therefore, we have
$$\mbox{meas}(Z(\phi(.;p),\eps)\cap]v-r,v+r[) \leq 
\overline{c}_i (\eps/\delta)^{1/i} \leq \overline{c} \, \eps^{1/k}/\delta$$
using Lemma \ref{Ls2}, where  $\overline{c}= \dis \max_{i=1,\cdots,k} \overline{c}_i$.
\\
 By compactness of  $K$, there exists a finite number 
 of such  sets $U_j=U(v_j,p_j)$ such that 
 $ \displaystyle{ K \subset \bigcup_{j=1}^{\nu} U_j}$. 
Thus, for each $p$, $Z(\phi(.;p),\eps)$ intersects at most $\nu$ 
 intervals $]v_j-r_j,v_j+r_j[$ where Lemma \ref{Ls2} is applied.  
 This allows to write
$\mbox{meas} (Z(\phi(.;p),\eps)) 
    \leq \nu c \, \eps^{1/k}/\delta $ for all $p$ and to conclude the proof.
$\CQFD$

\begin{Lemma}
 \label{Lpoly}
 Let $a(v)$ be the field $(v^1,v^2,\cdots,v^N)$ 
  then 
  $\alpha_{a(.)}= 1/N.$ 
\end{Lemma}
{\bf Proof.}
 From Lemma \ref{LsinfN}, we have yet $\alpha_{a(.)} \leq 1/N.$
  So, we just have to prove that $\alpha = 1/N$ satisfies (\ref{alpha})
 to conclude.
\\
 For all $v$, 
   $rank\{a'(v),\cdots,a^{(N)}(v)\} = N$, 
   thus it is impossible   to find  $\sigma  \in S^{N-1}$
 such that 
  $ \sigma \perp \{a'(v),\cdots,a^{(N)}(v)\}$.
 Let $\ph(v;u,\sigma)$ be $a(v) \cdot \sigma -u$.
 Since $\partial^{j}_v \ph(v;u,\sigma)=a^{(j)}(v) \cdot \sigma $ for $j \geq 1$,
 we have everywhere
 $
  \dis   \sum_{j=1}^N |\partial^{j}_v \ph(v;u,\sigma)| > 0.
$ 
\\
 Furthermore, for $|u| > 1 + a_{max}$,
  where  $a_{max}=\dis \sup_{|v|\leq A}|a(v)|$, 
  we have  $ |\ph(v;u,\sigma)| > 1$ 
for any $v \in [-A,A]$ and $\sigma \in S^{N-1}$.
Thus we can apply Lemma \ref{Lscompact} with $0 < \eps \leq 1$
 on the compact set 
 $[-A,A]_v\times [-a_{max}-1,a_{max} +1]_u \times S_\sigma^{N-1}$ which 
 concludes the proof with $\alpha_{a(.)}=1/N$.
$\CQFD$
\medskip\\  
{\bf Proof of Proposition \ref{maxalpha}.} 
 With Lemma \ref{LsinfN}, we have $\alpha_{opt}(N,1)\leq 1/N$.
 From Lemma \ref{Lpoly}, necessarily   $\alpha_{opt}(N,1)= 1/N$ which 
 concludes the proof.
 $\CQFD$

%%%%%%%%%%%%%%%%%%%%%%%%%%%%%%%%%%%%%%%%%%%%%%%%%%%%%%%%%%%%%%%%%%%%%%%
\Subsection{$M=N$}    
%%%%%%%%%%%%%%%%%%%%%%%%%%%%%%%%%%%%%%%%%%%%%%%%%%%%%%%%%%%%%%%%%%%%%%%
  The case when space dimension is equal  to  velocity dimension
 is the most physical one and then is very important. 
In this case, we can get the best smoothing effect with 
 $\alpha=1$. 
\begin{Prop}
 \label{alpha=1} For $N=M$, we have $\alpha_{opt}(N,N)= 1.$
\end{Prop}
{\bf Proof.} Since $\alpha \leq 1$, it is enough to find $a(.)$ 
  such that $\alpha=1$.
\\
 Let $a(.): \r^N_v \rightarrow \r^N_x$ be a global diffeomorphism, 
  $A>0$, $(u,\sigma) \in S^{N}$ and  
 $\ph(v)=a(v)\cdot \sigma -u$.
 Let $Z(\ph,\eps)=\{|v| \leq A, \, |\ph(v)|\leq \eps\}$. 
 Since $Da(v) \in GL_N(\r)$ and $\sigma \neq 0$, 
 then $\nabla_v \ph  \neq 0$ 
 and the set $Z(\ph,0)$ is empty or a manifold of dimension $N-1$. 
\\
Notice that for any $v$, there exists $(u,\sigma) \in S^N$ 
 such that $ a(v) \cdot \sigma -u = 0$, i.e. $Z(\ph,0)\neq \emptyset$.
  For instance, let $\widetilde{\sigma}$ belong to  $ S^{N-1}$
  and set $\widetilde{u}=a(v) \cdot\widetilde{\sigma} $,
  then  $(u,\sigma) = \dis \frac{1}{\sqrt{\widetilde{u}^2 +1}}
  ( \widetilde{u}, \widetilde{\sigma})$ satisfies the conditions.
\\
We thus consider that $Z(\ph,0)$ is not empty.
\\
There   exists $\delta$ such that 
 $ 0 < \delta <  |\nabla_v \ph(v)| < 1/\delta$
 for all $|v| \leq A, u ^2 +|\sigma|^2=1$. 
 \\
Using the mean inequality, we obtain  $\dis 
            \delta |v-v'| \leq |\ph(v) -\ph(v') | \leq \frac{|v-v'|}{\delta},$
 which implies for all $\eps < 1$, with 
$B(x,r) = \{y,\, |x-y|\leq r\} \subset \r^N$, that
\bqre 
  \bigcup_{z \in Z(\ph,0)} B(z,\delta \eps) 
   & \subset   Z(\ph,\eps )\subset& 
 \bigcup_{z \in Z(\ph,0)} B(z, \eps/\delta) 
\eqre
and
$Z(\ph,0)$ is diffeomorph to a piece of a hyperplane, so $meas (Z(\ph,\eps))$
 is of order $\eps$. More precisely, there exists a constant $C> 0$,
  only dependent on $A$,  $\delta$ and $||Da(.)||_{B(0,A)}$ 
 % and independent 
 such that $ 0 <C <\dis \frac{meas (Z(\ph,\eps))}{\eps} < C^{-1}.$
\\
Notice that  if $a(.)$ is a local diffeomorphism, $\alpha$ is still 
$1$.
$\CQFD$
\medskip \\
 Incidentally, we also  have $\alpha_{opt}(N,M)=1$ for all $M\geq N$.

%°°°°°°°°°°°°°°°°°°°°°°°°°°°°°°°°°°°°°°°°°°°°°°°°°°°°°°°°°°°°°°°°°°°°°°°°°°°°°
%°°°°°°°°°°°°°°°°°°°°°°°°°°°°°°°°°°°°°°°°°°°°°°°°°°°°°°°°°°°°°°°°°°°°°°°°°°°°°
%°°°°°°°°°°°°°°°°°°°°°°°°°°°°°°°°°°°°°°°°°°°°°°°°°°°°°°°°°°°°°°°°°°°°°°°°°°°°°
%°°°°°°°°°°°°°°°°°°°°°°°°°°°°°°°°°°°°°°°°°°°°°°°°°°°°°°°°°°°°°°°°°°°°°°°°°°°°°
%°°°°°°°°°°°°°°°°°°°°°°°°°°°°°°°°°°°°°°°°°°°°°°°°°°°°°°°°°°°°°°°°°°°°°°°°°°°°°
%°°°°°°°°°°°°°°°°°°°°°°°°°°°°°°°°°°°°°°°°°°°°°°°°°°°°°°°°°°°°°°°°°°°°°°°°°°°°°

\Section{Theorem in the $L^p$ framework}

Let us now deal with $L^p$ case. 
It will be an interpolation result of the $L^2$ obtained bound
and an estimate in $L^1$ using 
some operators in Hardy spaces.
We note ${\cal H}^1(\r^{N+1})$ the Hardy space and
${\cal H}^1(\r^{N}\times\r)$ the product Hardy space as done
in \cite{Bez} (see \cite{Ste} for more details about such spaces).

We will use the two following Propositions.
The first one is an interpolation result
(see \cite{Lin}, \cite{Bez} and \cite{BGP})
and the second one is about multiplier (\cite{Bez}).
\begin{Prop}[B\'ezard, Interpolation]
Let $T$ be a $\C$-linear operator, bounded in
$$ L^2(\r_t\times{\r^N}_x\times{\r^M}_v) \to W^{\beta,2}(\r_t\times\r^N_x),$$
and  in
$$L^1(\r^M_v, {\cal H}^1(\r^N \times \r)) \to {\cal H}^1(\r^{N+1}_{t,x}),$$
for some $\gamma \geq 0$.
Then $T$ is bounded 
$$L^p(\r_t\times\r^N_x\times\r^N_v) \to W^{s,p}(\r_t \times \r^N_x),$$
for $1<p\leq 2$, with $s=2\beta/{p'}$. 
\end{Prop}
\begin{Prop}[B\'ezard, Multiplier on ${\cal H}^1$]  \label{multiplicateur}
Let $m(y,y_{N+1})$ be a function of $(y,y_{n+1})\in \r^N\times\r$
which is $C^\infty$ out of
$[y=0 \textrm{ or } y_{N+1}=0]$, and verifying for all $\alpha$, $\beta$,
$$ |\partial_y^\alpha \partial_{y_{N+1}}^\beta m(y,y_{N+1}) |
\leq \frac{C_{\alpha \beta}}{|y|^\alpha |y_{N+1}|^\beta },$$
then $m$ defines a bounded Fourier multiplier on
${\cal H}^1 (\r^{N} \times \r)$.
\end{Prop}

\noindent {\bf Proof of Theorems 4 and 5.}\\
For Theorem 4 (respectively Theorem 5), we use 
the averaging lemma of Theorem 1 (respectively Theorem 2) which gives that
$T(f,g)=\rho_\psi$ is bounded from $L^2$ to 
$H^{\alpha/2}_{loc}$ (respectively $H^{1/\gamma}$).

We now focus on estimate in $L^1$.
We denote by ${\cal F}$ the Fourier transform with respect to $X$.
Taking this Fourier transform in $\ds b(v) \cdot \nabla_X f +F(X) \cdot \nabla_v f= g$,
we have
$${\cal F}(f)=\frac{{\cal F}(g)-
{\cal F}(F \cdot \nabla_vf)}{i(b(v)\cdot Y)}.$$
Let $\chi \in C^\infty_c(\r)$, $\chi(0)=1$, $\chi'(0)=0$ 
and $\chi''(0)\neq 0$ be an even,
non increasing function in $[0,+\infty[$.
We set $L$ such that $\textrm{supp} \chi \subset [-L,L]$.
We have
\begin{eqnarray*}
f(Y,v) &=& {\cal F}^{-1} \Big[ \chi(b(v)\cdot Y) {\cal F}(f)(Y,v)
+ (1-\chi(b(v)\cdot Y)) {\cal F}(f)(Y,v) \Big] \\
&=& {\cal F}^{-1} \Big[ \chi(b(v)\cdot Y) {\cal F}(f)(Y,v) \Big] \\
&&+ {\cal F}^{-1} \Big[ (1-\chi(b(v)\cdot Y)) \frac{{\cal F}(g)-
{\cal F}(F \cdot \nabla_vf)}{i(b(v)\cdot Y)} \Big],
\end{eqnarray*}
and then, in order to bound operator
$\ds f \mapsto \int_{\r^M} f(Y,v) \psi(v) \,dv$,
we have to bound the three following operators
\begin{equation} \label{dfnQ}
Q : f \mapsto \int_{\r^M} {\cal F}^{-1} 
\Big[ \chi(b(v)\cdot Y) {\cal F}(f)(Y,v) \Big] \psi(v) \,dv,
\end{equation}
\begin{equation} \label{dfnW}
W : g \mapsto \int_{\r^M} {\cal F}^{-1} 
\left[ \frac{1-\chi(b(v)\cdot Y)}{i(b(v)\cdot Y)}\,
{\cal F}(g)(Y,v) \right] \psi(v) \,dv
\end{equation}
and
\begin{equation} \label{dfnR}
R : f \mapsto -\int_{\r^M} {\cal F}^{-1} 
\left[ \frac{1-\chi(b(v)\cdot Y)}{i(b(v)\cdot Y)} \,
{\cal F}(F\cdot \nabla_v f)(Y,v) \right] \psi(v) \,dv.
\end{equation}
As in the classical case (by this we refer to \cite{Bez}, \cite{BGP}), we
transform the operators in order for them to involve only one
direction in $X$.
Indeed, the manipulation of product structure for Hardy space which depends on a moving
direction is difficult to deal with. Thus,
for any $v$, we take $R_v$ an orthogonal transform in
$\r^{N+1}$ such that $$R_v\left(\frac{b(v)}{|b(v)|}\right)=e_{N+1},$$ where $e_{N+1}$ is the very last 
vector of the canonical base,
and we set $$f_*(X,v)=f(R^{-1}_v(X),v)$$ and
$$Q_* f_*= Qf.$$ 
Since $f \mapsto f_*$ is an isometry on $L^p_{Xv}$,
we have now to study $Q_*$ instead of $Q$. We perform similar
transformations for the two other operators and we get $W_*$ and $R_*$.\\
For the two first operators, as in the classical proof, we have
$$\|Qf\|_{{\cal H}^1(\r^{N+1})} \leq C \|f\|_{L^1(\r^M_v, {\cal H}^1(\r^N\times
\r))},$$
and
$$\|Wg\|_{{\cal H}^1(\r^{N+1})} \leq C \|g\|_{L^1(\r^M_v, {\cal H}^1(\r^N\times
\r))}.$$
The new term is the third one (operator $R$).
We use the following rewrite of $R(f)$ in order to bound it. This is
\begin{eqnarray}
\nonumber
(Rf)(Y)&=&
-{\cal F}^{-1} \int_{\r^M} 
\left[ \frac{1-\chi(b(v)\cdot Y)}{i(b(v)\cdot Y)} \,
F\cdot \nabla_v {\cal F}(f)(Y,v) \right] \psi(v) \,dv \\
&=&{\cal F}^{-1}\left(  F \cdot \int_{\r^M} {\cal F}(f)(Y,v) \,
\nabla_v  \left[ \frac{1-\chi(b(v)\cdot Y)}{i(b(v)\cdot Y)} \,
\psi(v) \right]  \,dv \right) \nonumber \\
&=&{\cal F}^{-1}\left(  F \cdot \int_{\r^M} {\cal F}(f)(Y,v) \,
  \frac{1-\chi(b(v)\cdot Y)}{i(b(v)\cdot Y)} \,
\nabla_v \psi(v)   \,dv \right) \nonumber \\
&&+{\cal F}^{-1}\left(  F \cdot \int_{\r^M} {\cal F}(f)(Y,v) \,
m_0(b(v)\cdot Y) \nabla_v (b(v) \cdot Y) \,
\psi(v)  \,dv \right) \nonumber \\
\label{reecriture}
\end{eqnarray}
with 
\begin{equation} \label{dfnmo}
m_0(y)= \frac{-y\chi'(y)-1+\chi(y)}{iy^2}.
\end{equation}
We denote by ${\cal F}(R_1f)$ and ${\cal F}(R_2f)$ 
the two terms of this decomposition. We perform as previously
orthogonal transformations and we have to study the obtained $(R_1)_*$ and $(R_2)_*$.\\
The term $(R_1)_*$ is the same than $W_*$
but with $\nabla_v \psi$ instead of $\psi$. Thus we have the same result
thanks to the regularity assumption on $\psi$.\\
Now, setting $T=m_0 \nabla_v$,  we have
\begin{eqnarray*}
(R_2)_*(f_*)(Y) &=&
 F \cdot \int_{\r^M}  {\cal F}^{-1}\bigg( {\cal F}(f_*)(R_v(Y),v) \,
T\Big(b(v)\cdot Y\Big)  \bigg) \,
\psi(v)  \,dv \\
 &=&
 F \cdot \int_{\r^M}  {\cal F}^{-1}\bigg( {\cal F}(f_*)(R_v(Y),v) \,
T\Big(R_v(b(v))\cdot R_v(Y)\Big)  \bigg) \,
\psi(v)  \,dv \\
&=&
 F \cdot \int_{\r^M}  {\cal F}^{-1}\bigg( {\cal F}(f_*)(R_v(Y),v) \,
T\Big(|b(v)| e_{N+1}\cdot R_v(Y)\Big)  \bigg) \,
\psi(v)  \,dv, 
\end{eqnarray*}
thus, setting $T_j=m_0 \partial_{v_j}$,  we get
\begin{eqnarray*}
&&\|(R_2)_*(f_*)\|_{{\cal H}^1(\r^{N+1})} \\
&\leq &  \sum_j |F_j| \int_{\r^M} \left\| {\cal F}^{-1} 
\bigg( {\cal F}(f_*)(R_v(Y),v) \,
T_j\Big(|b(v)| e_{N+1}\cdot R_v(Y)\Big)  \bigg)
\right\|_{{\cal H}^1(\r^{N+1})} |\psi(v)| \,dv\\
&\leq &  \sum_j |F_j| \int_{\r^M} \left\| {\cal F}^{-1} 
\bigg( {\cal F}(f_*)(Y,v) \,
T_j\Big(|b(v)| e_{N+1}\cdot Y\Big)  \bigg)
\right\|_{{\cal H}^1(\r^{N+1})} |\psi(v)| \,dv\\
&\leq &  C_1 \sum_j |F_j| \int_{\r^M} \left\| {\cal F}^{-1} 
\bigg( {\cal F}(f_*)(Y,v) \,
T_j\Big(|b(v)| e_{N+1}\cdot Y\Big)  \bigg)
 \right\|_{{\cal H}^1(\r^{N}\times\r)} |\psi(v)| \,dv,
\end{eqnarray*}
using the invariance under orthogonal transformation in ${\cal H}^1(\r^{N+1})$
and thanks to the continuous injection of
${\cal H}^1(\r^{N}\times\r)$
in ${\cal H}^1(\r^{N+1})$.\\
We use now 
Proposition \ref{multiplicateur}
 with the term
$$m_j(y,y_{N+1})=
T_j(|b(v)| e_{N+1} \cdot Y)=
m_0(|b(v)| y_{N+1}) \partial_{v_j} (|b(v)|) y_{N+1}, \quad \textrm{ for }
j=1,\cdots, M.$$
Those terms rewrite
$$m_j(y,y_{N+1})=m_0(|b(v)| y_{N+1}) \frac{a(v) \cdot \partial_{v_j} a(v)}{|b(v)|} y_{N+1}.$$
Now $m_0(z) \tod_{z \to 0} -\frac{1}{2i} \chi''(0),$
therefore $m_0$ is $C^\infty$.
The terms in (\ref{dfnmo}) with $\chi$ have a compact support and the other term is $1/y^2$, then
every derivatives of $m_0$ is bounded at infinity.\\
\noindent We differentiate  $m_j$ with respect to $y_{N+1}$, it gives
\begin{eqnarray*}
 \partial_{y_{N+1}}^{k} m_j(y,y_{N+1}) & =& 
\frac{a(v) \cdot \partial_{v_j} a(v)}{|b(v)|}
\Big( m_0^{(k)} ( |b(v)| y_{N+1}) 
|b(v)|^k y_{N+1}\\
&& \hphantom{\frac{a(v) \cdot \partial_{v_j} a(v)}{|b(v)|} aa}  + k m_0^{(k-1)} ( |b(v)| y_{N+1}) 
|b(v)|^{k-1} \Big).
\end{eqnarray*}
There exists  some constants $C$ and $C_k$ such that
$$|b(v)| \leq C, \qquad |b(v)|^{k-2} |a(v) \cdot \partial_{v_j} a(v)| \leq C_k$$
for $v$ in the compact support of $\psi$.
Thus 
$$\left|\partial_{y_{N+1}}^{k} m_j(y,y_{N+1}) \right| |y_{N+1}|^k
\leq C_k \left(  C m_0^{(k)} ( |b(v)| y_{N+1}) 
 y_{N+1} + k m_0^{(k-1)} ( |b(v)| y_{N+1})  \right).$$
For $|y_{N+1}| \geq (R+1)/C$, we have $m_0^{(j)} ( |b(v)| y_{N+1})=0$ for any $j$,
and then $m_0^{(k)} ( |b(v)| y_{N+1})  y_{N+1} + k m_0^{(k-1)} ( |b(v)| y_{N+1})
= 0$ for $|y_{N+1}| \geq (R+1)/C$.\\
Furthermore $\ds |m_0^{(k)} ( |b(v)| y_{N+1})  y_{N+1} + k m_0^{(k-1)} ( |b(v)| y_{N+1})|
\leq \|m_0^{(k)}\|_\infty \frac{R+1}{C} + k \|m_0^{(k-1)}\|_\infty$
for $|y_{N+1}| < (R+1)/C$.
Finally, for any $(y,y_{N+1})$, we get
$$\left|\partial_{y_{N+1}}^{k} m_j(y,y_{N+1}) \right| |y_{N+1}|^k
\leq C_k \left(   \|m_0^{(k)}\|_\infty (R+1) + k \|m_0^{(k-1)}\|_\infty
 \right)$$
uniformly with respect to $v$ in the support of $\psi$.
Then, we can apply Proposition \ref{multiplicateur} to get the boundary of  $(R_2)_*$.\\
The interpolation result concludes, since $\beta=\alpha/2$
(respectively $\beta=1/\gamma$), that
the obtained regularity is $s=\alpha/p'$
(respectively $s=2/(\gamma p')$).\CQFD

\bigskip

{\bf Acknowledgments.} We thank Gilles Lebeau and Jeffrey Rauch for fruitful
discussions on oscillatory integrals.
We also thank referees for there useful comments.

%{\bf R\'ef\'erences }

\end{document}